\documentclass{amsart}

\usepackage{graphicx}
\usepackage{mathtools}

\usepackage{tikz}

\usetikzlibrary{arrows.meta,positioning}

\DeclarePairedDelimiter\floor{\lfloor}{\rfloor}

 \newtheorem{theorem}{Theorem}
\newtheorem{proposition}[theorem]{Proposition}
 
 \newtheorem{remark}[theorem]{Remark}

\title{Planar tropical caustics:\\trivalency and convexity}

\author{Mikhail Shkolnikov}

\email{m.shkolnikov@math.bas.bg}

\thanks{Institute of Mathematics and Informatics, Bulgarian Academy of Sciences, Sofia.\\
Supported by the Simons Foundation, grant SFI-MPS-T-Institutes-00007697, and\\ the Ministry of Education and Science of the Republic of Bulgaria, grant DO1-239/10.12.2024.}

\newcommand{\K}{\mathcal{K}}

\begin{document}

\begin{abstract}
Tropical caustic of a convex domain on the plane is a canonically associated tropical analytic curve inside the domain. In this note we give a graphical proof for the classification of its intermediate vertices, implying in particular that they are always trivalent. Apart from that we explain how various known examples of tropical caustics are constructed and discuss the possibility of relaxing the convexity condition for the domain.
\end{abstract}

\maketitle

\section{Introduction}
\begin{figure}[h!]
\includegraphics[width=0.9\textwidth]{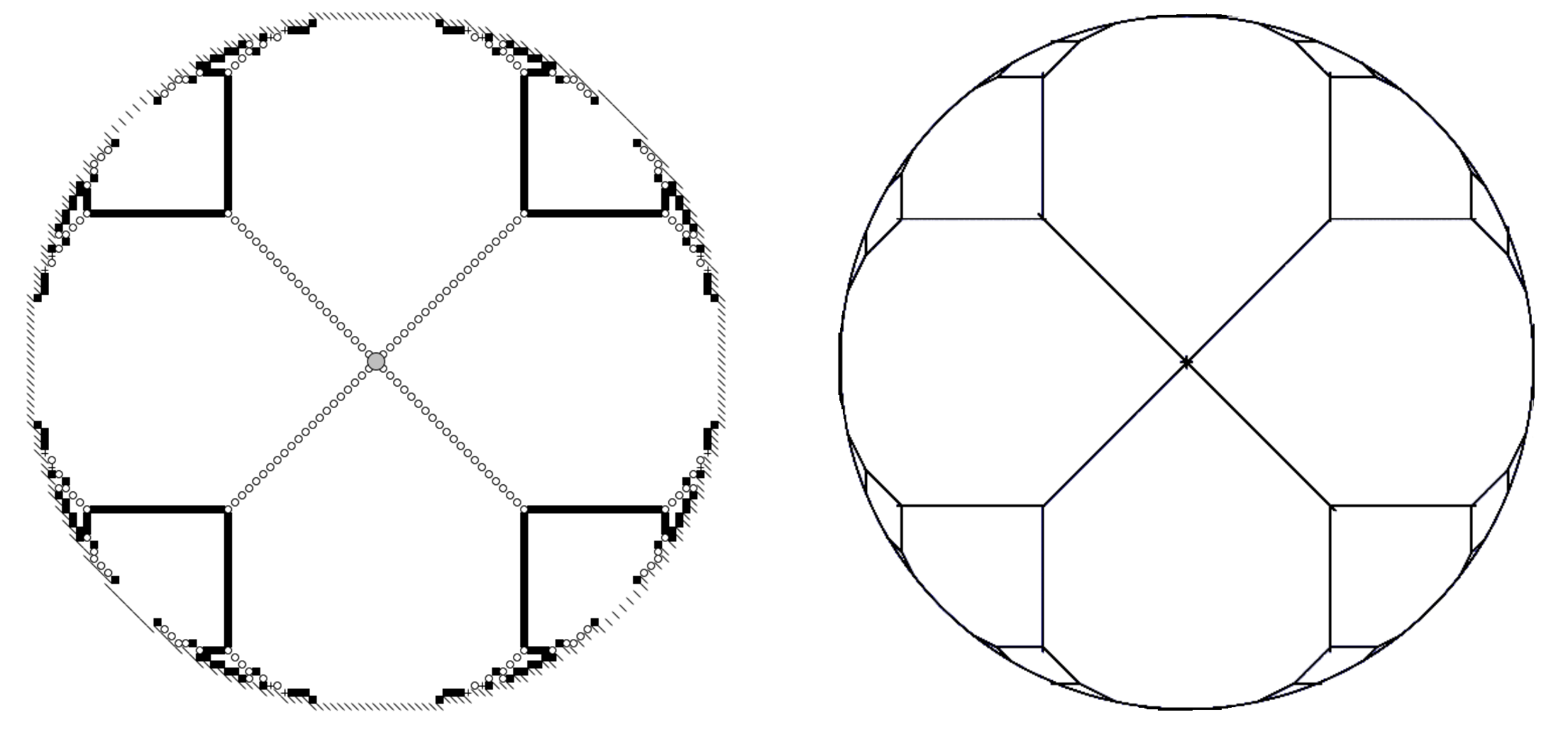}
\caption{Right: the result of dropping a single grain to the center of a disc on the maximal stable state of the sandpile model. Left: tropical caustic of the disc.}
\label{fig_causticdisc}
\end{figure}

For a convex domain $\Phi\subset\mathbb{R}^2,$ denote by $\K_\Phi\subset\Phi$ its tropical caustic. Several constructions of $\K_\Phi$ will be reviewed in this note. If $\Phi$ is compact, $\K_\Phi\subset\Phi$ is a rational tropical analytic curve, i.e. it is a possibly infinite tree with rational slope edges enhanced with weights such that the balancing condition at every vertex is satisfied. This tree has infinitely many edges unless $\Phi$ is a polygonal domain with rational slope sides. The first caustic of a domain $\Phi$ with a smooth boundary ever observed was the tropical caustic of a disc (see Figure \ref{fig_causticdisc}).  This was done in the context of proving certain scaling limit theorems for the sandpile model \cite{KS15,KS16}, giving rise to the continuous version of the model introduced in \cite{KL18} and further investigated in \cite{KP23,S23,V25}. Tropical curves arising in this scaling limit minimize action (the integral of the tropical series) and symplectic area (in case of lattice polygons). The latter functional is a weighted sum of lengths of edges of a tropical curve, making it into a solution of a Steiner-type problem \cite{S17}. Recently, it was shown in \cite{CP} that branching points for solutions of planar Gilbert-Steiner problem (under quite mild assumptions) are trivalent, and even more recently a similar fact under some conditions but in higher dimension was established in \cite{C25}. Although, it is not clear how exactly tropical caustics are related to these minimizers, there are many parallels between the theories. One of them is that intermediate vertices of tropical caustics are trivalent. This fact for unimodular polygon $\Phi$ (proven in \cite{KS18,S17}) is instrumental in establishing the tropical sandpile scaling limit theorems. For a general convex domain $\Phi,$ it was proven in \cite{MS23} by means of tropical trigonometry as a part of systematic treatment of planar tropical wave fronts and their caustics. In this note, we give an alternative graphical proof of the intermediate trivalency phenomenon (see Theorem \ref{thm_nonfinal}). The parallel with optimal branched transport theory, suggests that the convexity restriction on $\Phi$ could be potentially dropped, this possibility is discussed in the last section.

\begin{figure}

\centering
\begin{tikzpicture}[
    >=Stealth,
    node/.style={draw, thick, align=center, inner xsep=10pt, inner ysep=8pt}
]
\node[node] (DS) at (-3, 2) {Tropical distance series\\[2pt]
  $\mathcal{F}_{\Phi}:\Phi\rightarrow[0,+\infty)$};
\node[node] (WF) at ( 3, 1) {Tropical wave front\\[2pt]
  $t\geq 0\mapsto\partial\Phi(t)$};
\node[node] (CA) at ( 0,-2) {Tropical caustic\\[2pt]
  $\mathcal{K}_{\Phi}\subset\Phi$};

\draw[<->, thick] (WF) -- (DS);
\draw[<->, thick] (WF) -- (CA);
\draw[<->, thick] (DS) -- (CA);

\node[align=center] at (0,0) {\it Tropical \\\it Optics};
\end{tikzpicture}

\caption{The three fundamental entities of tropical optics. Each one uniquely restores the other two.}
\end{figure}
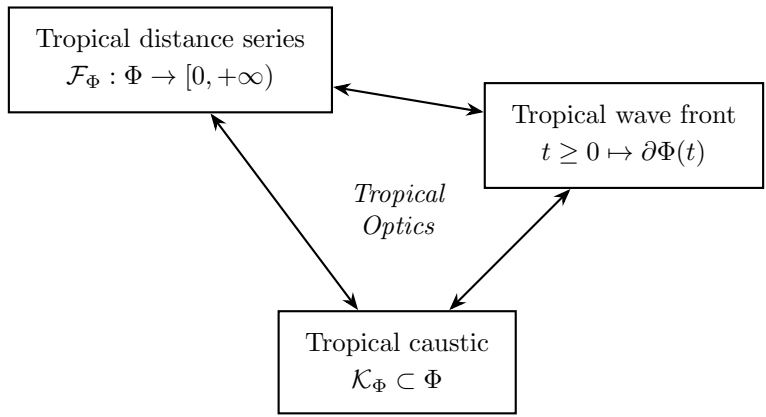

\section{Definitions and properties}
The definitions and results listed in this section follow \cite{MS23}. For simplicity, we start with the case of compact domains and mention what happens in the non-compact case later on.
\subsection{Compact domains}

Let $N$ be a rank two lattice spanning a plane $N_\mathbb{R}=N\otimes\mathbb{R}.$ Denote by $M$ the dual lattice $\operatorname{Hom}(N,\mathbb{Z}).$ A half-plane $H\subset N_\mathbb{R}$ is said to have a rational slope if it can be presented as $$H=\{p\in N_\mathbb{R}:\lambda(p)\geq c\}$$ for some $\lambda\in M$ and $c\in\mathbb{R}.$ Note that we can assume that $\lambda$ is primitive. In this case, a propagation by $t\in\mathbb{R}$ of $H$ is given by  $$H(t)=\{p\in N_\mathbb{R}:\lambda(p)\geq c+t\}.$$

Let $\Phi\subset N_\mathbb{R}$ be a compact convex domain. A half-plane $H$ is said to be a support half-plane of $\Phi$ if $\Phi\subset H$ and $\partial H\cap \partial \Phi\neq \emptyset.$ Denote by $\mathcal{H}_\Phi$ the set of all rational slope support half-planes of $\Phi.$ The tropical propagation of $\Phi$ in time $t\in\mathbb{R}$ is defined as $$\Phi(t)=\bigcap\limits_{H\in\mathcal{H}_\Phi} H(t).$$
The {\it tropical wave front} of $\partial\Phi$ at time $t$ is $\partial \Phi(t).$

\begin{remark}For all $t\in\mathbb{R},$ $\Phi(t)$ is convex. 
\end{remark}

\begin{proposition}For $t\leq 0,$ $\Phi(t)=\Phi.$
\end{proposition}

\begin{proposition}For $t>0,$ $\Phi(t)\subset\Phi^\circ.$
\end{proposition}

\begin{theorem}[Huygens' principle] For all $t,s>0,$ $$(\Phi(t))(s)=\Phi(t+s).$$
\end{theorem}

\begin{proposition}
For all non-empty $\Phi$ there exist $t_\Phi\geq 0$ such that $\Phi(t)$ is not empty, but has an empty interior. 
\end{proposition}

We call such $t_\Phi$ {\it the final time} of $\Phi.$ 

\begin{proposition}For $t>0,$ $\Phi(t)$ is a polygon with rational slope sides. \label{prop_polygon}
\end{proposition}

This enables us to employ terminology from toric geometry through a generalized version of Delzant's theorem (such as in \cite{LT97}). Namely, we may interpret $\Phi(t)$ as a moment polygon of some possibly singular toric surface $S(t)$ with symplectic form $\omega(t)$ for $t\in (0,t_\Phi).$ We call a time $t\in(0,t_\Phi)$ a critical time if for all $\epsilon\in(0,t)$ the dual fans of $\Phi(t)$ and $\Phi(t-\epsilon)$ are different. In particular, we can assume that between critical times the surface $S(t)$ doesn't change. 

\begin{proposition}
The set of critical times of $\Phi$ has at most one accumulation point. If exists, this point is $0.$
\end{proposition}

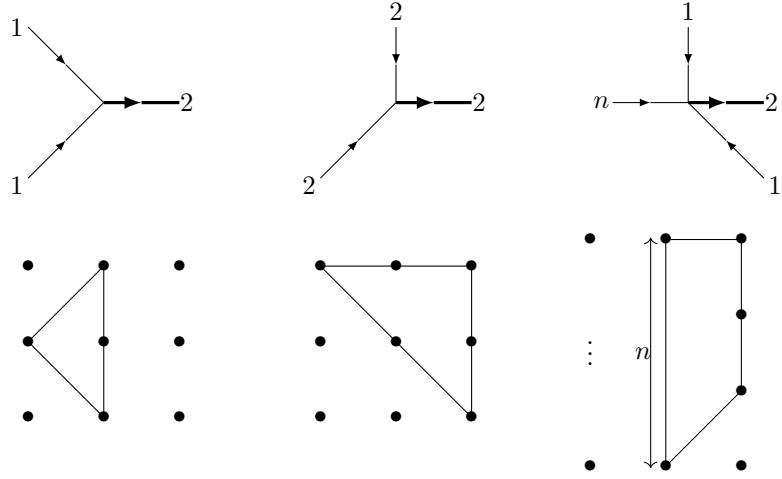
\begin{figure}[h!]
\begin{tikzpicture}
\node at (-1.15,1) {$1$};
\draw[-latex](-1,1)--(-0.5,0.5);
\draw(-0.5,0.5)--(0,0);

\node at (-1.15,-1.1) {$1$};
\draw[-latex](-1,-1)--(-0.5,-0.5);
\draw(-0.5,-0.5)--(0,0);

\node at (1.1,0) {$2$};
\draw[-latex,very thick](0,0)--(0.5,0);
\draw[very thick](0.5,0)--(1,0);

\begin{scope}[yshift=-90]
\draw(-1,0)--(0,1)--(0,-1)--(-1,0);
\node at (-1,1) {$\bullet$};
\node at (-1,0) {$\bullet$};
\node at (-1,-1) {$\bullet$};
\node at (0,1) {$\bullet$};
\node at (0,0) {$\bullet$};
\node at (0,-1) {$\bullet$};
\node at (1,1) {$\bullet$};
\node at (1,0) {$\bullet$};
\node at (1,-1) {$\bullet$};
\end{scope}

\begin{scope}[xshift=110]
\node at (0,1.2) {$2$};
\draw[-latex](0,1)--(0,0.5);
\draw(0,0.5)--(0,0);

\node at (-1.15,-1.1) {$2$};
\draw[-latex](-1,-1)--(-0.5,-0.5);
\draw(-0.5,-0.5)--(0,0);

\node at (1.1,0) {$2$};
\draw[-latex,very thick](0,0)--(0.5,0);
\draw[very thick](0.5,0)--(1,0);

\begin{scope}[yshift=-90]
\draw(-1,1)--(1,1)--(1,-1)--(-1,1);

\node at (-1,1) {$\bullet$};
\node at (-1,0) {$\bullet$};
\node at (-1,-1) {$\bullet$};
\node at (0,1) {$\bullet$};
\node at (0,0) {$\bullet$};
\node at (0,-1) {$\bullet$};
\node at (1,1) {$\bullet$};
\node at (1,0) {$\bullet$};
\node at (1,-1) {$\bullet$};
\end{scope}

\end{scope}

\begin{scope}[xshift=220]
\node at (0,1.2) {$1$};
\draw[-latex](0,1)--(0,0.5);
\draw(0,0.5)--(0,0);

\node at (-1.15,0) {$n$};
\draw[-latex](-1,0)--(-0.5,0);
\draw(-0.5,0)--(0,0);

\node at (1.15,-1.1) {$1$};
\draw[-latex](1,-1)--(0.5,-0.5);
\draw(0.5,-0.5)--(0,0);

\node at (1.1,0) {$2$};
\draw[-latex,very thick](0,0)--(0.5,0);
\draw[very thick](0.5,0)--(1,0);

\begin{scope}[yshift=-80,xshift=20]
\draw(0,-1)--(-1,-2)--(-1,1)--(0,1)--(0,-1);

\draw[<->](-1.2,1.025)--(-1.2,-2.025);
\node at (-1.3,-0.5) {$n$};

\node at (-2,-0.4) {$\vdots$};

\node at (-1,1) {$\bullet$};
\node at (-1,-2) {$\bullet$};

\node at (-2,1) {$\bullet$};
\node at (-2,-2) {$\bullet$};

\node at (0,1) {$\bullet$};
\node at (0,0) {$\bullet$};
\node at (0,-1) {$\bullet$};
\node at (0,-2) {$\bullet$};

\end{scope}

\end{scope}

\end{tikzpicture}
\caption{Possible schemes of collisions in the particle process (defined in Section \ref{sec_trival}) at the final time in the case when a mass two particle is formed, its trajectory is traced in bold. The configuration on the right is realizable for all natural $n.$ Below: dual polygons of corresponding vertices where collisions take place. \\ Note that, in this case, there is a second collision with another weight two particle formed at the same (final) time. These two particles move towards each other along the final locus $\Phi(t_\Phi)$. Each pair of the schemes above can occur independently at the two sides of this edge, if the domain $\Phi$ is compact. For a non-compact domain with parallel rational asymptotes of the boundary) the final locus is a ray, i.e. the caustic final singularity type is determined by a single picture out of these types.}\label{fig_finaledge}
\end{figure}

\begin{theorem}\label{thm_An}
For any compact $\Phi$ and $t\in(0,t_\Phi)$ the toric surface $S(t)$ may have only $A_n$-type singularities. 
For a critical time $t$ and small $\varepsilon>0$, the transformation from $S(t-\varepsilon)$ to $S(t)$ is a collection of contractions (blowing down) of isolated boundary divisors $D_j\subset S(t-\varepsilon)$.
Each exceptional divisor $D_j$ contains at most one singular point of $S(t-\varepsilon)$.
The point resulting in contraction of $D_j$ is a non-singular point of $S(t)$.
\end{theorem}

Since $S(t)$ for $t\in(0,t_\Phi)$ has only mild singularities, its canonical class $K_{S(t)}$ is well defined as a class in the second cohomology of $S(t)$ with integer coefficients. Here is how the de Rham cohomology class of the symplectic form changes with time.

\begin{theorem}Between critical times,
\begin{equation} \label{eq_canonicalevolution}{\frac{d}{dt}}[\omega(t)]=2\pi\cdot K_{S(t)}\in H^2_{dR}(S(t)).\tag{$\star$}\end{equation}
\label{thm_canev} 
\end{theorem}
 
A point of the wave front is said to be special if there is more than one support half-plane in it. In view of Proposition \ref{prop_polygon}, for $t\in(0,t_\Phi)$ only vertices of $\Phi(t)$ are special and special points of $\Phi(0)=\Phi$ are points where its boundary is not smooth. Note that all points of $\Phi(t_\Phi)$ are special. Therefore, we have a definition of a tropical caustic $\mathcal{K}_\Phi$ of $\Phi$ as a set --- it is a locus of all special points of all $\Phi(t).$

\begin{theorem}$\mathcal{K}_\Phi$ is a tropical analytic curve, i.e. it is a graph with straight edges of rational slopes, balanced at every vertex in the interior of $\Phi$ with weights given by the following rule: an $A_n$-vertex of $\Phi(t)$ traces a weight $n+1$ edge of $\mathcal{K}_\Phi$ and if $\Phi(t_\Phi)$ is a segment, it has weight $2$.
\end{theorem}

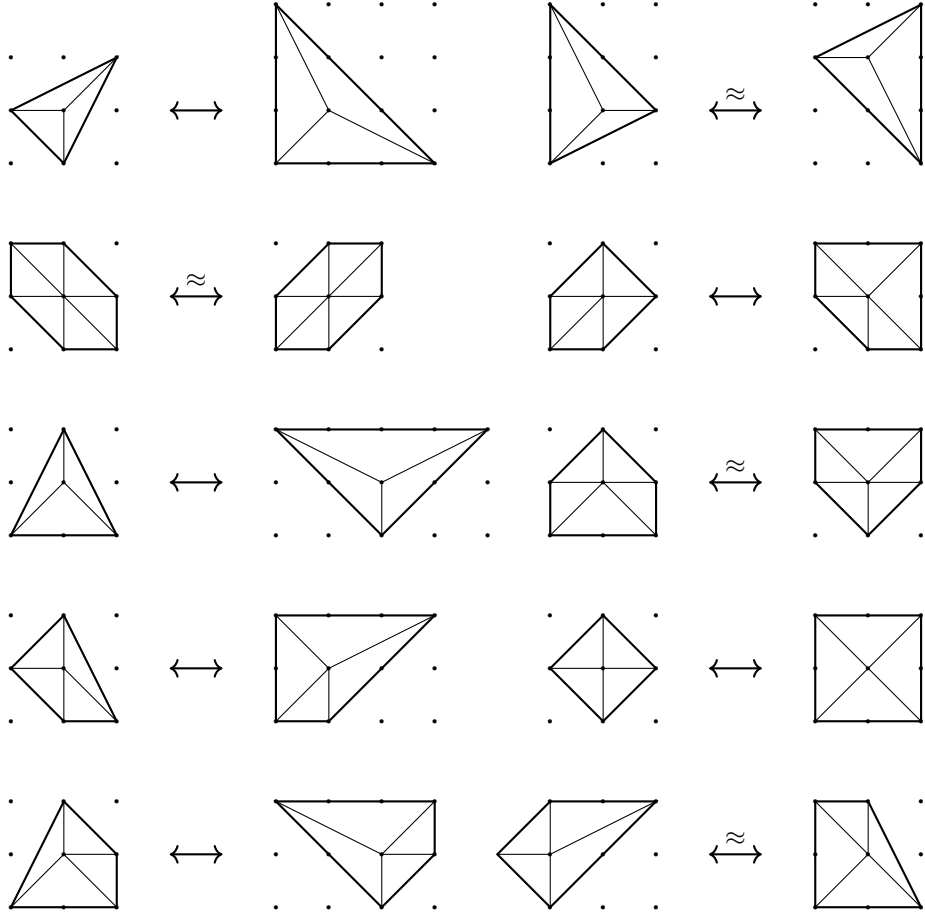
\begin{figure}[h!]
\begin{tikzpicture}

\begin{scope}[scale=0.7]

\node at (0,0) {$\boldsymbol{\cdot}$};
\node at (1,0) {$\boldsymbol{\cdot}$};
\node at (2,0) {$\boldsymbol{\cdot}$};
\node at (0,1) {$\boldsymbol{\cdot}$};
\node at (1,1) {$\boldsymbol{\cdot}$};
\node at (2,1) {$\boldsymbol{\cdot}$};
\node at (0,2) {$\boldsymbol{\cdot}$};
\node at (1,2) {$\boldsymbol{\cdot}$};
\node at (2,2) {$\boldsymbol{\cdot}$};

\draw[thick](1,0)--(2,2)--(0,1)--cycle;

\draw(1,1)--(1,0);
\draw(1,1)--(2,2);
\draw(1,1)--(0,1);

\draw[thick,<->] (3,1)--(4,1);

\begin{scope}[xshift=142.5]
\node at (0,0) {$\boldsymbol{\cdot}$};
\node at (1,0) {$\boldsymbol{\cdot}$};
\node at (2,0) {$\boldsymbol{\cdot}$};
\node at (3,0) {$\boldsymbol{\cdot}$};
\node at (0,1) {$\boldsymbol{\cdot}$};
\node at (1,1) {$\boldsymbol{\cdot}$};
\node at (2,1) {$\boldsymbol{\cdot}$};
\node at (3,1) {$\boldsymbol{\cdot}$};
\node at (0,2) {$\boldsymbol{\cdot}$};
\node at (1,2) {$\boldsymbol{\cdot}$};
\node at (2,2) {$\boldsymbol{\cdot}$};
\node at (3,2) {$\boldsymbol{\cdot}$};
\node at (0,3) {$\boldsymbol{\cdot}$};
\node at (1,3) {$\boldsymbol{\cdot}$};
\node at (2,3) {$\boldsymbol{\cdot}$};
\node at (3,3) {$\boldsymbol{\cdot}$};

\draw[thick](0,0)--(0,3)--(3,0)--cycle;

\draw(1,1)--(0,0);
\draw(1,1)--(3,0);
\draw(1,1)--(0,3);
\end{scope} 

\begin{scope}[xshift=290]
\node at (0,0) {$\boldsymbol{\cdot}$};
\node at (1,0) {$\boldsymbol{\cdot}$};
\node at (2,0) {$\boldsymbol{\cdot}$};
\node at (0,1) {$\boldsymbol{\cdot}$};
\node at (1,1) {$\boldsymbol{\cdot}$};
\node at (2,1) {$\boldsymbol{\cdot}$};
\node at (0,2) {$\boldsymbol{\cdot}$};
\node at (1,2) {$\boldsymbol{\cdot}$};
\node at (2,2) {$\boldsymbol{\cdot}$};
\node at (0,3) {$\boldsymbol{\cdot}$};
\node at (1,3) {$\boldsymbol{\cdot}$};
\node at (2,3) {$\boldsymbol{\cdot}$};

\draw[thick](0,0)--(0,3)--(2,1)--cycle;

\draw(1,1)--(0,0);
\draw(1,1)--(0,3);
\draw(1,1)--(2,1);

\draw[thick,<->] (3,1)--(4,1);

\node at (3.5,1.3) {$\approx$};

\begin{scope}[xshift=142.5]
\node at (0,0) {$\boldsymbol{\cdot}$};
\node at (1,0) {$\boldsymbol{\cdot}$};
\node at (2,0) {$\boldsymbol{\cdot}$};
\node at (0,1) {$\boldsymbol{\cdot}$};
\node at (1,1) {$\boldsymbol{\cdot}$};
\node at (2,1) {$\boldsymbol{\cdot}$};
\node at (0,2) {$\boldsymbol{\cdot}$};
\node at (1,2) {$\boldsymbol{\cdot}$};
\node at (2,2) {$\boldsymbol{\cdot}$};
\node at (0,3) {$\boldsymbol{\cdot}$};
\node at (1,3) {$\boldsymbol{\cdot}$};
\node at (2,3) {$\boldsymbol{\cdot}$};

\draw[thick](2,0)--(2,3)--(0,2)--cycle;

\draw(1,2)--(2,0);
\draw(1,2)--(2,3);
\draw(1,2)--(0,2);
\end{scope}

\end{scope}

\begin{scope}[yshift=-100]

\node at (0,0) {$\boldsymbol{\cdot}$};
\node at (1,0) {$\boldsymbol{\cdot}$};
\node at (2,0) {$\boldsymbol{\cdot}$};
\node at (0,1) {$\boldsymbol{\cdot}$};
\node at (1,1) {$\boldsymbol{\cdot}$};
\node at (2,1) {$\boldsymbol{\cdot}$};
\node at (0,2) {$\boldsymbol{\cdot}$};
\node at (1,2) {$\boldsymbol{\cdot}$};
\node at (2,2) {$\boldsymbol{\cdot}$};

\draw[thick](0,1)--(0,2)--(1,2)--(2,1)--(2,0)--(1,0)--cycle;

\draw(1,1)--(0,1);
\draw(1,1)--(0,2);
\draw(1,1)--(1,2);
\draw(1,1)--(2,1);
\draw(1,1)--(2,0);
\draw(1,1)--(1,0);

\draw[thick,<->] (3,1)--(4,1);

\node at (3.5,1.3) {$\approx$};

\begin{scope}[xshift=142.5]
\node at (0,0) {$\boldsymbol{\cdot}$};
\node at (1,0) {$\boldsymbol{\cdot}$};
\node at (2,0) {$\boldsymbol{\cdot}$};
\node at (0,1) {$\boldsymbol{\cdot}$};
\node at (1,1) {$\boldsymbol{\cdot}$};
\node at (2,1) {$\boldsymbol{\cdot}$};
\node at (0,2) {$\boldsymbol{\cdot}$};
\node at (1,2) {$\boldsymbol{\cdot}$};
\node at (2,2) {$\boldsymbol{\cdot}$};

\draw[thick](0,0)--(0,1)--(1,2)--(2,2)--(2,1)--(1,0)--cycle;

\draw(1,1)--(0,0);
\draw(1,1)--(0,1);
\draw(1,1)--(1,2);
\draw(1,1)--(2,2);
\draw(1,1)--(2,1);
\draw(1,1)--(1,0);

\end{scope} 

\begin{scope}[xshift=290]
\node at (0,0) {$\boldsymbol{\cdot}$};
\node at (1,0) {$\boldsymbol{\cdot}$};
\node at (2,0) {$\boldsymbol{\cdot}$};
\node at (0,1) {$\boldsymbol{\cdot}$};
\node at (1,1) {$\boldsymbol{\cdot}$};
\node at (2,1) {$\boldsymbol{\cdot}$};
\node at (0,2) {$\boldsymbol{\cdot}$};
\node at (1,2) {$\boldsymbol{\cdot}$};
\node at (2,2) {$\boldsymbol{\cdot}$};

\draw[thick](0,0)--(0,1)--(1,2)--(2,1)--(1,0)--cycle;

\draw(1,1)--(0,0);
\draw(1,1)--(0,1);
\draw(1,1)--(1,2);
\draw(1,1)--(2,1);
\draw(1,1)--(1,0);

\draw[thick,<->] (3,1)--(4,1);

\begin{scope}[xshift=142.5]
\node at (0,0) {$\boldsymbol{\cdot}$};
\node at (1,0) {$\boldsymbol{\cdot}$};
\node at (2,0) {$\boldsymbol{\cdot}$};
\node at (0,1) {$\boldsymbol{\cdot}$};
\node at (1,1) {$\boldsymbol{\cdot}$};
\node at (2,1) {$\boldsymbol{\cdot}$};
\node at (0,2) {$\boldsymbol{\cdot}$};
\node at (1,2) {$\boldsymbol{\cdot}$};
\node at (2,2) {$\boldsymbol{\cdot}$};

\draw[thick](0,1)--(0,2)--(2,2)--(2,0)--(1,0)--cycle;

\draw(1,1)--(0,1);
\draw(1,1)--(0,2);
\draw(1,1)--(2,2);
\draw(1,1)--(2,0);
\draw(1,1)--(1,0);

\end{scope}

\end{scope}

\end{scope}

\begin{scope}[yshift=-200]

\node at (0,0) {$\boldsymbol{\cdot}$};
\node at (1,0) {$\boldsymbol{\cdot}$};
\node at (2,0) {$\boldsymbol{\cdot}$};
\node at (0,1) {$\boldsymbol{\cdot}$};
\node at (1,1) {$\boldsymbol{\cdot}$};
\node at (2,1) {$\boldsymbol{\cdot}$};
\node at (0,2) {$\boldsymbol{\cdot}$};
\node at (1,2) {$\boldsymbol{\cdot}$};
\node at (2,2) {$\boldsymbol{\cdot}$};

\draw[thick](0,0)--(1,2)--(2,0)--cycle;

\draw(1,1)--(0,0);;
\draw(1,1)--(1,2);
\draw(1,1)--(2,0);

\draw[thick,<->] (3,1)--(4,1);

\begin{scope}[xshift=142.5]
\node at (0,0) {$\boldsymbol{\cdot}$};
\node at (1,0) {$\boldsymbol{\cdot}$};
\node at (2,0) {$\boldsymbol{\cdot}$};
\node at (3,0) {$\boldsymbol{\cdot}$};
\node at (4,0) {$\boldsymbol{\cdot}$};
\node at (0,1) {$\boldsymbol{\cdot}$};
\node at (1,1) {$\boldsymbol{\cdot}$};
\node at (2,1) {$\boldsymbol{\cdot}$};
\node at (3,1) {$\boldsymbol{\cdot}$};
\node at (4,1) {$\boldsymbol{\cdot}$};
\node at (0,2) {$\boldsymbol{\cdot}$};
\node at (1,2) {$\boldsymbol{\cdot}$};
\node at (2,2) {$\boldsymbol{\cdot}$};
\node at (3,2) {$\boldsymbol{\cdot}$};
\node at (4,2) {$\boldsymbol{\cdot}$};

\draw[thick](0,2)--(4,2)--(2,0)--cycle;

\draw(2,1)--(0,2);
\draw(2,1)--(4,2);
\draw(2,1)--(2,0);

\end{scope} 

\begin{scope}[xshift=290]
\node at (0,0) {$\boldsymbol{\cdot}$};
\node at (1,0) {$\boldsymbol{\cdot}$};
\node at (2,0) {$\boldsymbol{\cdot}$};
\node at (0,1) {$\boldsymbol{\cdot}$};
\node at (1,1) {$\boldsymbol{\cdot}$};
\node at (2,1) {$\boldsymbol{\cdot}$};
\node at (0,2) {$\boldsymbol{\cdot}$};
\node at (1,2) {$\boldsymbol{\cdot}$};
\node at (2,2) {$\boldsymbol{\cdot}$};

\draw[thick](0,0)--(0,1)--(1,2)--(2,1)--(2,0)--cycle;

\draw(1,1)--(0,0);
\draw(1,1)--(0,1);
\draw(1,1)--(1,2);
\draw(1,1)--(2,1);
\draw(1,1)--(2,0);

\draw[thick,<->] (3,1)--(4,1);

\node at (3.5,1.3) {$\approx$};

\begin{scope}[xshift=142.5]
\node at (0,0) {$\boldsymbol{\cdot}$};
\node at (1,0) {$\boldsymbol{\cdot}$};
\node at (2,0) {$\boldsymbol{\cdot}$};
\node at (0,1) {$\boldsymbol{\cdot}$};
\node at (1,1) {$\boldsymbol{\cdot}$};
\node at (2,1) {$\boldsymbol{\cdot}$};
\node at (0,2) {$\boldsymbol{\cdot}$};
\node at (1,2) {$\boldsymbol{\cdot}$};
\node at (2,2) {$\boldsymbol{\cdot}$};

\draw[thick](0,1)--(0,2)--(2,2)--(2,1)--(1,0)--cycle;

\draw(1,1)--(0,1);
\draw(1,1)--(0,2);
\draw(1,1)--(2,2);
\draw(1,1)--(2,1);
\draw(1,1)--(1,0);

\end{scope}

\end{scope}

\end{scope}

\begin{scope}[yshift=-300]

\node at (0,0) {$\boldsymbol{\cdot}$};
\node at (1,0) {$\boldsymbol{\cdot}$};
\node at (2,0) {$\boldsymbol{\cdot}$};
\node at (0,1) {$\boldsymbol{\cdot}$};
\node at (1,1) {$\boldsymbol{\cdot}$};
\node at (2,1) {$\boldsymbol{\cdot}$};
\node at (0,2) {$\boldsymbol{\cdot}$};
\node at (1,2) {$\boldsymbol{\cdot}$};
\node at (2,2) {$\boldsymbol{\cdot}$};

\draw[thick](0,1)--(1,2)--(2,0)--(1,0)--cycle;

\draw(1,1)--(0,1);
\draw(1,1)--(1,2);
\draw(1,1)--(2,0);
\draw(1,1)--(1,0);

\draw[thick,<->] (3,1)--(4,1);

\begin{scope}[xshift=142.5]
\node at (0,0) {$\boldsymbol{\cdot}$};
\node at (1,0) {$\boldsymbol{\cdot}$};
\node at (2,0) {$\boldsymbol{\cdot}$};
\node at (3,0) {$\boldsymbol{\cdot}$};
\node at (0,1) {$\boldsymbol{\cdot}$};
\node at (1,1) {$\boldsymbol{\cdot}$};
\node at (2,1) {$\boldsymbol{\cdot}$};
\node at (3,1) {$\boldsymbol{\cdot}$};
\node at (0,2) {$\boldsymbol{\cdot}$};
\node at (1,2) {$\boldsymbol{\cdot}$};
\node at (2,2) {$\boldsymbol{\cdot}$};
\node at (3,2) {$\boldsymbol{\cdot}$};

\draw[thick](0,0)--(0,2)--(3,2)--(1,0)--cycle;

\draw(1,1)--(0,0);
\draw(1,1)--(0,2);
\draw(1,1)--(3,2);
\draw(1,1)--(1,0);

\end{scope} 

\begin{scope}[xshift=290]
\node at (0,0) {$\boldsymbol{\cdot}$};
\node at (1,0) {$\boldsymbol{\cdot}$};
\node at (2,0) {$\boldsymbol{\cdot}$};
\node at (0,1) {$\boldsymbol{\cdot}$};
\node at (1,1) {$\boldsymbol{\cdot}$};
\node at (2,1) {$\boldsymbol{\cdot}$};
\node at (0,2) {$\boldsymbol{\cdot}$};
\node at (1,2) {$\boldsymbol{\cdot}$};
\node at (2,2) {$\boldsymbol{\cdot}$};

\draw[thick](0,1)--(1,2)--(2,1)--(1,0)--cycle;

\draw(1,1)--(1,0);
\draw(1,1)--(0,1);
\draw(1,1)--(2,1);
\draw(1,1)--(1,2);

\draw[thick,<->] (3,1)--(4,1);

\begin{scope}[xshift=142.5]
\node at (0,0) {$\boldsymbol{\cdot}$};
\node at (1,0) {$\boldsymbol{\cdot}$};
\node at (2,0) {$\boldsymbol{\cdot}$};
\node at (0,1) {$\boldsymbol{\cdot}$};
\node at (1,1) {$\boldsymbol{\cdot}$};
\node at (2,1) {$\boldsymbol{\cdot}$};
\node at (0,2) {$\boldsymbol{\cdot}$};
\node at (1,2) {$\boldsymbol{\cdot}$};
\node at (2,2) {$\boldsymbol{\cdot}$};

\draw[thick](0,0)--(0,2)--(2,2)--(2,0)--cycle;

\draw(1,1)--(0,0);
\draw(1,1)--(2,0);
\draw(1,1)--(0,2);
\draw(1,1)--(2,2);
\end{scope}

\end{scope}

\end{scope}

\begin{scope}[yshift=-400]

\node at (0,0) {$\boldsymbol{\cdot}$};
\node at (1,0) {$\boldsymbol{\cdot}$};
\node at (2,0) {$\boldsymbol{\cdot}$};
\node at (0,1) {$\boldsymbol{\cdot}$};
\node at (1,1) {$\boldsymbol{\cdot}$};
\node at (2,1) {$\boldsymbol{\cdot}$};
\node at (0,2) {$\boldsymbol{\cdot}$};
\node at (1,2) {$\boldsymbol{\cdot}$};
\node at (2,2) {$\boldsymbol{\cdot}$};

\draw[thick](0,0)--(1,2)--(2,1)--(2,0)--cycle;

\draw(1,1)--(0,0);;
\draw(1,1)--(1,2);
\draw(1,1)--(2,1);
\draw(1,1)--(2,0);

\draw[thick,<->] (3,1)--(4,1);

\begin{scope}[xshift=142.5]
\node at (0,0) {$\boldsymbol{\cdot}$};
\node at (1,0) {$\boldsymbol{\cdot}$};
\node at (2,0) {$\boldsymbol{\cdot}$};
\node at (3,0) {$\boldsymbol{\cdot}$};
\node at (0,1) {$\boldsymbol{\cdot}$};
\node at (1,1) {$\boldsymbol{\cdot}$};
\node at (2,1) {$\boldsymbol{\cdot}$};
\node at (3,1) {$\boldsymbol{\cdot}$};
\node at (0,2) {$\boldsymbol{\cdot}$};
\node at (1,2) {$\boldsymbol{\cdot}$};
\node at (2,2) {$\boldsymbol{\cdot}$};
\node at (3,2) {$\boldsymbol{\cdot}$};

\draw[thick](0,2)--(3,2)--(3,1)--(2,0)--cycle;

\draw(2,1)--(0,2);
\draw(2,1)--(3,2);
\draw(2,1)--(3,1);
\draw(2,1)--(2,0);

\end{scope} 

\begin{scope}[xshift=290]
\node at (0,0) {$\boldsymbol{\cdot}$};
\node at (1,0) {$\boldsymbol{\cdot}$};
\node at (2,0) {$\boldsymbol{\cdot}$};
\node at (0,1) {$\boldsymbol{\cdot}$};
\node at (1,1) {$\boldsymbol{\cdot}$};
\node at (2,1) {$\boldsymbol{\cdot}$};
\node at (0,2) {$\boldsymbol{\cdot}$};
\node at (1,2) {$\boldsymbol{\cdot}$};
\node at (2,2) {$\boldsymbol{\cdot}$};

\draw[thick](-1,1)--(0,2)--(2,2)--(0,0)--cycle;

\draw(0,1)--(-1,1);
\draw(0,1)--(0,2);
\draw(0,1)--(2,2);
\draw(0,1)--(0,0);

\draw[thick,<->] (3,1)--(4,1);

\node at (3.5,1.3) {$\approx$};

\begin{scope}[xshift=142.5]
\node at (0,0) {$\boldsymbol{\cdot}$};
\node at (1,0) {$\boldsymbol{\cdot}$};
\node at (2,0) {$\boldsymbol{\cdot}$};
\node at (0,1) {$\boldsymbol{\cdot}$};
\node at (1,1) {$\boldsymbol{\cdot}$};
\node at (2,1) {$\boldsymbol{\cdot}$};
\node at (0,2) {$\boldsymbol{\cdot}$};
\node at (1,2) {$\boldsymbol{\cdot}$};
\node at (2,2) {$\boldsymbol{\cdot}$};

\draw[thick](0,0)--(0,2)--(1,2)--(2,0)--cycle;

\draw(1,1)--(0,0);
\draw(1,1)--(0,2);
\draw(1,1)--(1,2);
\draw(1,1)--(2,0);

\end{scope}

\end{scope}

\end{scope}

\end{scope}
\end{tikzpicture}
\caption{Sixteen types of lattice polygons (up to automorphisms of the lattice) with a single lattice point in the interior together with their tropical caustics. Each caustic consists of segments connecting the central point to the vertices. The polygons are paired by duality, the weight of an edge of the caustic is the length of the corresponding side of the dual polygon. Note that some of these types are self-dual, in a sense of identifying the lattices $M$ and $N.$ Observe that the sums of lattice perimeters of two dual polygons are always equal to 12 and that the self-dual types are exactly the 3-,4-,5- and 6-gon with the perimeter 6. This is an instance of a general Noether-type formula established in \cite{MS23} equating the sum of the perimeter of the domain boundary and its caustic to the length of the final segment (multiplied by $4$), plus  the final time $t_\Phi$ (multiplied by $12$). \\ These sixteen types classify possible final singularities of tropical caustics in the case when the final locus $\Phi(t_\Phi)$ of the wave front propagation is a point.
}
\label{fig_reflexivepoly}
\end{figure}

Singularities of tropical caustics fall into three categories:
\begin{itemize}
\item Initial, i.e. at $t=0$ -- described by lattice points in dual cones;
\item Intermediate, at critical $t$ -- a single type for every natural $n$ (Figure \ref{fig_nonfinal});
\item Final, i.e. at $t=t_\Phi$ -- in the case of $\Phi(t_\Phi)$ equal to a point, coincides with the famous list of $16$ reflexive polygons (Figure \ref{fig_reflexivepoly}), and in the case of $\Phi(t_\Phi)$ equal to a segment, has one sporadic case and a single type for every natural $n=0,1,2,\dots$ (Figure \ref{fig_finaledge}).
\end{itemize}

\begin{theorem}\label{thm_causticcoordinates}
The intermediate and final singularity types of $\mathcal{K}_\Phi,$ the lengths of its edges and their cyclic order restore $\Phi$ up to translations and automorphisms of the lattice $N$.
\end{theorem}

In addition, we have the following suplementary structural statement.

\begin{remark}\label{rem_crittimes_and_lengths}
For a vertex $p$ of the tropical caustic $\mathcal{K}_\Phi,$ the correponding critical time is equal to $\mathcal{F}_\Phi(p).$ For a pair of adjacent vertices $p_1$ and $p_2$, the length of the caustic edge connecting them is equal to the difference between their critical times.
\end{remark}

This result implies that one can extract invariants of domains (such as area or perimeter) from the moduli of the caustic, which are equivalent to the multi-set of critical times. The first instance of such a formula was for the area of the disc \cite{KS17}, which boils down to the following identity
$$2\sum (\sqrt{a^2+b^2}+\sqrt{c^2+d^2}-\sqrt{(a+c)^2+(b+d)^2})^2=4-\pi,$$
where the sum is taken over all $a,b,c,d\in\mathbb{Z}_{\geq 0}$ satisfying $ad-bc=1.$ Here, $\pi$ is the area of the unit disc, and $4$ is the area of the square in which this disc is inscribed, this square should be thought as a ``symplectic minimal model'' (the definition of this concept is given in the Examples section) of the disc. The summands correspond to vertices of the caustic, which in turn are in one-to-one correspondence with a sequence of corner cuts (symplectic blow-ups, in toric terminology). The values in the summation, are squares of critical times.  Formulas for other domains and other invariants were obtained in \cite{KS19} and even more formulas appear in \cite{K24,K25}.

\subsection{Non-compact domains}
In the case of $\Phi$ non-compact convex domain, one big difference with the compact case is that the final time might be infinite. Also, it may happen that the domain has no support half-planes of rational slope. This situation occurs only in three cases: when $\Phi$ is the whole plane, when it is a half-plane of irrational slope and when it is a strip of irrational slope. In such cases, it is reasonable to pose the tropical wave front to be constant. One motivation for such a choice is that there is a way to get the continuous wave front propagation as a the scaling limit of the discrete one (a step of which consists of taking the convex hull of the set of lattice points in the interior of the domain).

Even the case of half-planes of rational slope is particular since the wave front has no special points, and thus the tropical caustic is empty. In all other situations, the caustic is non-empty, however, in contrast with the compact case it can be disconnected. Such phenomena arise when one desingularizes (an example of a resolution of singularities is shown in Figure \ref{fig_resolution}) a non-unimodular cone -- such a caustic will be a union of disjoint rays.

Caustics of cones are fundamental, but their structure is not yet fully understood. One reason why they are important is that tropical caustic of a polygon can be reconstructed from tropical caustics of the corners of this polygon via the particle process described in the next section. Another reason is due to the relation to the quantum toric planes and their moduli space. Currently, we understand the structure of this space only in the case when one of the sides of the cone has rational slope. In this situtation, it is parametrized by a single quantity in $\mathbb{R}\slash\mathbb{Z}$ called the tropical cotangent. In \cite{MS23}, continued fraction expressions for tropical cotangent in terms of tropical caustic and wave front were obtained, showing in particular that any sequence of caustic weights and wave front gradients are realizable. It is not clear, how this result is extended to the case when both sides of the cone have irrational slope.

\begin{figure}
\includegraphics[width=0.95\textwidth]{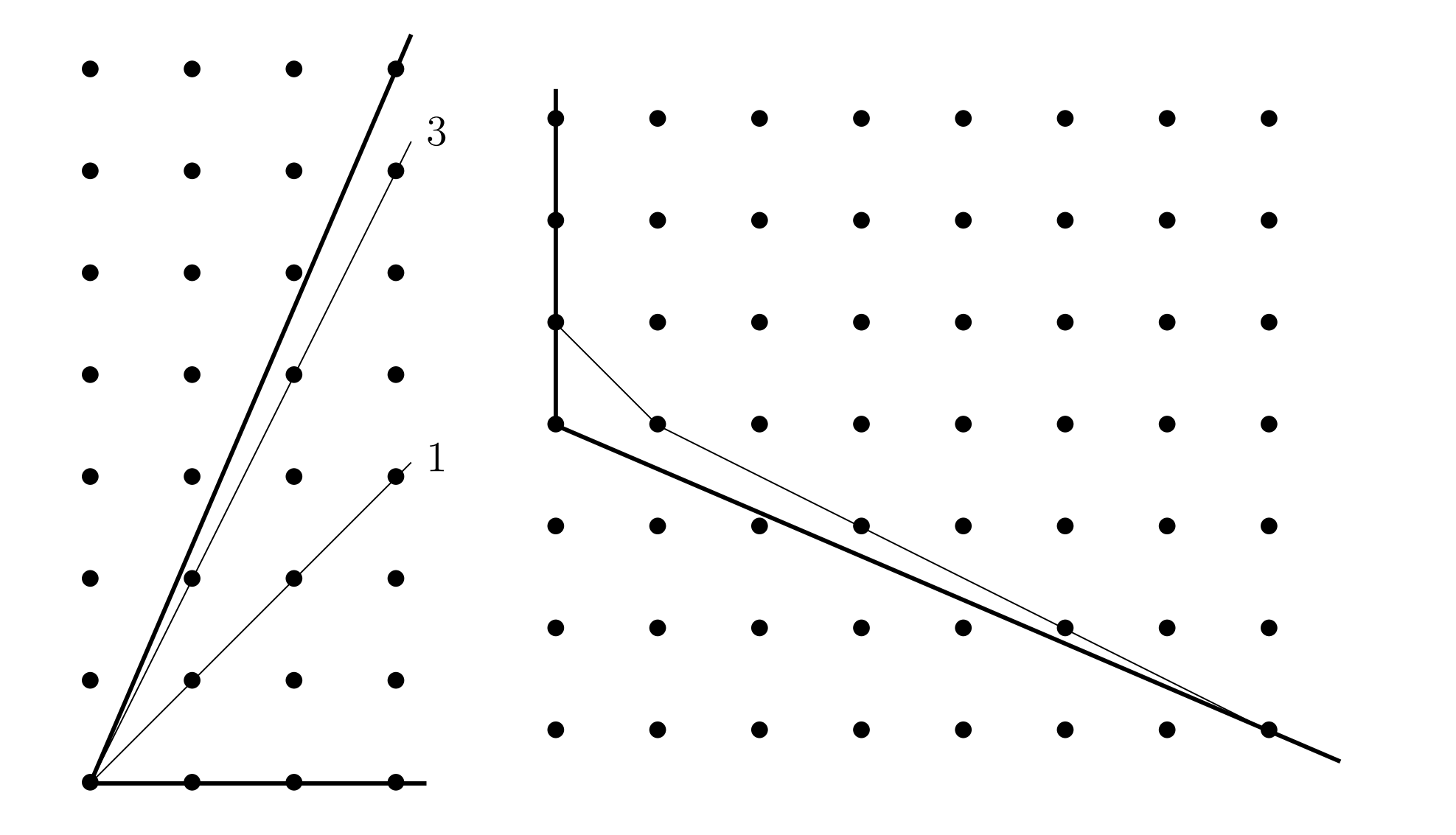}
\caption{The rule for constructing a caustic of a cone. One should take the convex hull of the set of non-zero lattice points in the dual cone. The finite edges in the resulting polygonal domain are in one-to-one correspondence (via orthogonality) with the rays in the caustic of the original cone, with the length of a segment being equal to the weight of the ray.}
\label{fig_conerule}
\end{figure}

\section{Trivalency at intermediate vertices}
\label{sec_trival}
Now we proceed with the main theorem and its (mostly graphical) proof. First, we need two more interpretations of the caustic. The first interpretation is that it is the corner locus of $\mathcal{F}_\Phi,$ a particular $\Phi$-tropical series \cite{KS18}, which is equal to the minimum of all support functions of $\Phi$ with integral gradient -- thus we may call it the ``tropical distance series'' of $\Phi$. The term ``corner locus'' stands simply for all points where $\mathcal{F}_\Phi$ is not linear, i.e. two or more support functions equal to the minimum. Another way stating this relation is to say that the caustic curve $\mathcal{K}_\Phi$ is defined by ``$\mathcal{F}_\Phi=0$'' in the hyperfield terminology. It is important to observe that the wave front evolution can also be expressed in terms of the tropical series, i.e. $\partial\Phi(t)=\mathcal{F}_\Phi^{-1}(t).$ 

Another way to construct the caustic of a polygonal domain is via the total trajectory of the particle process, where particles represent vertices of the wave front. Start with sending particles from the corners of the polygonal domain according to the rule of Figure \ref{fig_conerule}, each particle going along the caustic ray of the corner with primitive velocity and having the mass equal to the weight of the corresponding ray. When two or more particles meet, either they annihilate each other if the total momentum\footnote{By momentum of a praticle here we mean the product of its mass, i.e. the weight of the caustic segment, and the primitive velocity vector in the direction of the ray.} is zero, or a new particle is formed to have the conservation of the momentum, again going with a primitive velocity. The theorem below tells that at a non-final collision (i.e. when it happened before the final time $t_\Phi$ of the wave front propagation), only two particles can meet and the new particle has mass $1.$ 
   
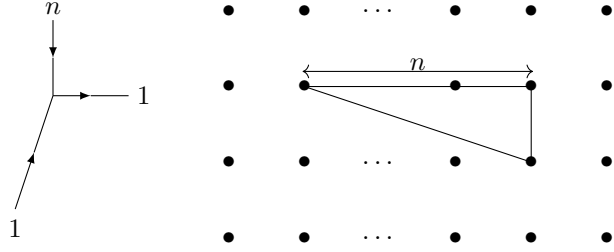
\begin{figure}[h!]
\begin{tikzpicture}
\node at (0,1.15) {$n$};
\draw[-latex](0,1)--(0,0.5);
\draw(0,0.5)--(0,0);

\node at (-0.5,-1.75) {$1$};
\draw[-latex](-0.5,-1.5)--(-0.25,-0.75);
\draw(-0.25,-0.75)--(0,0);

\node at (1.2,0) {$1$};
\draw[-latex](0,0)--(0.5,0);
\draw(0.5,0)--(1,0);

\begin{scope}[xshift=180,yshift=-25]
\draw(0,0)--(0,1)--(-3,1)--(0,0);

\draw[<->](0.025,1.2)--(-3.025,1.2);
\node at (-1.5,1.3) {$n$};

\node at (1,-1) {$\bullet$};
\node at (1,0) {$\bullet$};
\node at (1,1) {$\bullet$};
\node at (1,2) {$\bullet$};

\node at (0,-1) {$\bullet$};
\node at (0,0) {$\bullet$};
\node at (0,1) {$\bullet$};
\node at (0,2) {$\bullet$};

\node at (-1,-1) {$\bullet$};
\node at (-1,0) {$\bullet$};
\node at (-1,1) {$\bullet$};
\node at (-1,2) {$\bullet$};

\node at (-3,-1) {$\bullet$};
\node at (-3,0) {$\bullet$};
\node at (-3,1) {$\bullet$};
\node at (-3,2) {$\bullet$};

\node at (-2,-1) {$\dots$};
\node at (-2,-0) {$\dots$};

\node at (-2,2) {$\dots$};

\node at (-4,-1) {$\bullet$};
\node at (-4,0) {$\bullet$};
\node at (-4,1) {$\bullet$};
\node at (-4,2) {$\bullet$};
\end{scope}
\end{tikzpicture}
\caption{Illustration to Theorem \ref{thm_nonfinal}. The only possible scheme of a collision at a non-final time, for every natural $n,$ involving a single particle of weight one and a single particle of weight $n$ and resulting in a single particle of weight one. On the right: the dual $A_n$-type triangle of the corresponding trivalent vertex, unique up to automorphisms of the lattice.}\label{fig_nonfinal}
\end{figure}

\begin{theorem}\label{thm_nonfinal}
There exist a unique, up to automorphisms of the lattice, local model of a non-final collision (shown on Figure \ref{fig_nonfinal}). It involves exactly two particles, only one of them can have mass greater than one.
\end{theorem}

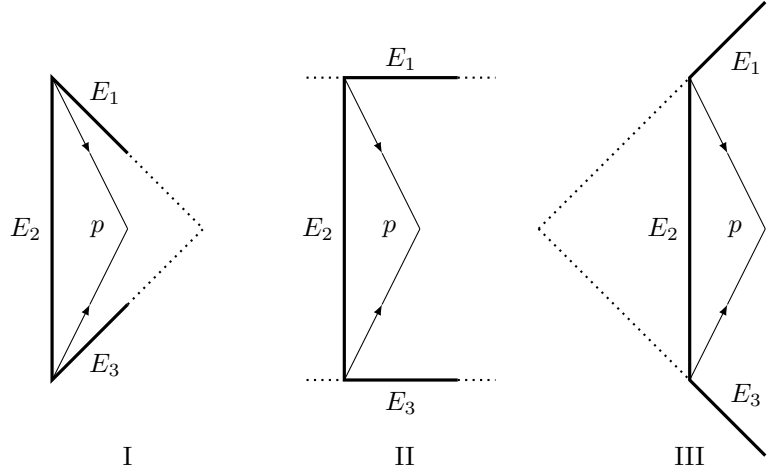
\begin{figure}
\begin{tikzpicture}

\draw[-latex](0,-2)--(0.5,-1);
\draw(0.5,-1)--(1,0);
\draw[-latex](0,2)--(0.5,1);
\draw(0.5,1)--(1,0);
\node at (0.6,0) {$p$};

\draw[very thick](1,-1)--(0,-2)--(0,2)--(1,1);
\draw[dotted, thick](1,-1)--(2,0)--(1,1);

\node at (0.7,1.8) {$E_1$};
\node at (-0.35,0) {$E_2$};
\node at (0.7,-1.8) {$E_3$};

\node at (1,-3){I};

\begin{scope}[xshift=110]

\draw[-latex](0,-2)--(0.5,-1);
\draw(0.5,-1)--(1,0);
\draw[-latex](0,2)--(0.5,1);
\draw(0.5,1)--(1,0);
\node at (0.6,0) {$p$};

\draw[very thick](1.5,-2)--(0,-2)--(0,2)--(1.5,2);
\draw[dotted, thick](-0.5,-2)--(0,-2);
\draw[dotted, thick](1.5,-2)--(2,-2);
\draw[dotted, thick](-0.5,2)--(0,2);
\draw[dotted, thick](1.5,2)--(2,2);

\node at (0.75,2.25) {$E_1$};
\node at (-0.35,0) {$E_2$};
\node at (0.75,-2.3) {$E_3$};

\node at (0.8,-3){II};
\end{scope}

\begin{scope}[xshift=240]
\draw[-latex](0,-2)--(0.5,-1);
\draw(0.5,-1)--(1,0);
\draw[-latex](0,2)--(0.5,1);
\draw(0.5,1)--(1,0);
\node at (0.6,0) {$p$};

\draw[very thick](1,-3)--(0,-2)--(0,2)--(1,3);
\draw[dotted, thick](-2,0)--(0,-2);
\draw[dotted, thick](-2,0)--(0,2);

\node at (0.75,2.2) {$E_1$};
\node at (-0.35,0) {$E_2$};
\node at (0.75,-2.2) {$E_3$};

\node at (0,-3){III};
\end{scope}

\end{tikzpicture}
\caption{Illustration to the proof of Theorem \ref{thm_nonfinal}. Three distinct geometric configurations of a first collision according to where the lines extending sides $E_1$ and $E_3$ intersect with respect to the line extending the side $E_2.$ Note that more than two particles may collide at point $p,$ but only two trajectories of consecutive particles are shown.}\label{fig_firstcollision}
\end{figure}

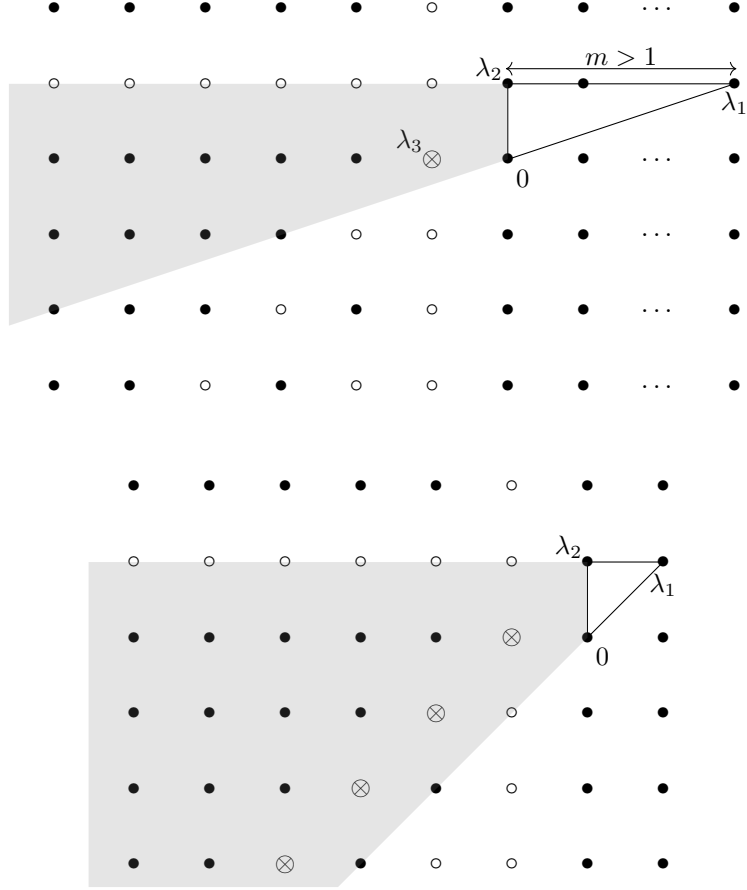
\begin{figure}
\begin{tikzpicture}
\draw(0,0)--(0,1)--(3,1)--(0,0);

\draw[<->](-0.025,1.2)--(3.025,1.2);
\node at (1.5,1.35) {$m>1$};

\node at (0.2,-0.25) {$0$};
\node at (3,0.75) {$\lambda_1$};
\node at (-0.25,1.2) {$\lambda_2$};
\node at (-1.3,0.25) {$\lambda_3$};

\node at (0,-3) {$\bullet$};
\node at (0,-2) {$\bullet$};
\node at (0,-1) {$\bullet$};
\node at (0,0) {$\bullet$};
\node at (0,1) {$\bullet$};
\node at (0,2) {$\bullet$};

\node at (1,-3) {$\bullet$};
\node at (1,-2) {$\bullet$};
\node at (1,-1) {$\bullet$};
\node at (1,0) {$\bullet$};
\node at (1,1) {$\bullet$};
\node at (1,2) {$\bullet$};

\node at (2,-3) {$\dots$};
\node at (2,-2) {$\dots$};
\node at (2,-1) {$\dots$};
\node at (2,-0) {$\dots$};
\node at (2,2) {$\dots$};

\node at (3,-3) {$\bullet$};
\node at (3,-2) {$\bullet$};
\node at (3,-1) {$\bullet$};
\node at (3,0) {$\bullet$};
\node at (3,1) {$\bullet$};
\node at (3,2) {$\bullet$};

\node at (-1,-3) {$\circ$};
\node at (-1,-2) {$\circ$};
\node at (-1,-1) {$\circ$};
\node at (-1,0) {$\otimes$};
\node at (-1,1) {$\circ$};
\node at (-1,2) {$\circ$};

\node at (-2,-3) {$\circ$};
\node at (-2,-2) {$\bullet$};
\node at (-2,-1) {$\circ$};
\node at (-2,0) {$\bullet$};
\node at (-2,1) {$\circ$};
\node at (-2,2) {$\bullet$};

\node at (-3,2) {$\bullet$};
\node at (-3,1) {$\circ$};
\node at (-3,0) {$\bullet$};
\node at (-3,-1) {$\bullet$};
\node at (-3,-2) {$\circ$};
\node at (-3,-3) {$\bullet$};

\node at (-4,-3) {$\circ$};
\node at (-4,-2) {$\bullet$};
\node at (-4,-1) {$\bullet$};
\node at (-4,0) {$\bullet$};
\node at (-4,1) {$\circ$};
\node at (-4,2) {$\bullet$};

\node at (-5,-3) {$\bullet$};
\node at (-5,-2) {$\bullet$};
\node at (-5,-1) {$\bullet$};
\node at (-5,0) {$\bullet$};
\node at (-5,1) {$\circ$};
\node at (-5,2) {$\bullet$};

\node at (-6,-3) {$\bullet$};
\node at (-6,-2) {$\bullet$};
\node at (-6,-1) {$\bullet$};
\node at (-6,0) {$\bullet$};
\node at (-6,1) {$\circ$};
\node at (-6,2) {$\bullet$};

\draw[fill=gray, draw=none, opacity=0.2](0,0)--(-6.6,-2.2)--(-6.6,1)--(0,1)--(0,0);

\begin{scope}[xshift=30,yshift=-180]

\draw(0,0)--(0,1)--(1,1)--(0,0);

\node at (0.2,-0.25) {$0$};
\node at (1,0.7) {$\lambda_1$};
\node at (-0.25,1.2) {$\lambda_2$};

\node at (0,-3) {$\bullet$};
\node at (0,-2) {$\bullet$};
\node at (0,-1) {$\bullet$};
\node at (0,0) {$\bullet$};
\node at (0,1) {$\bullet$};
\node at (0,2) {$\bullet$};

\node at (1,-3) {$\bullet$};
\node at (1,-2) {$\bullet$};
\node at (1,-1) {$\bullet$};
\node at (1,0) {$\bullet$};
\node at (1,1) {$\bullet$};
\node at (1,2) {$\bullet$};

\node at (-1,-3) {$\circ$};
\node at (-1,-2) {$\circ$};
\node at (-1,-1) {$\circ$};
\node at (-1,0) {$\otimes$};
\node at (-1,1) {$\circ$};
\node at (-1,2) {$\circ$};

\node at (-2,-3) {$\circ$};
\node at (-2,-2) {$\bullet$};
\node at (-2,-1) {$\otimes$};
\node at (-2,0) {$\bullet$};
\node at (-2,1) {$\circ$};
\node at (-2,2) {$\bullet$};

\node at (-3,2) {$\bullet$};
\node at (-3,1) {$\circ$};
\node at (-3,0) {$\bullet$};
\node at (-3,-1) {$\bullet$};
\node at (-3,-2) {$\otimes$};
\node at (-3,-3) {$\bullet$};

\node at (-4,-3) {$\otimes$};
\node at (-4,-2) {$\bullet$};
\node at (-4,-1) {$\bullet$};
\node at (-4,0) {$\bullet$};
\node at (-4,1) {$\circ$};
\node at (-4,2) {$\bullet$};

\node at (-5,-3) {$\bullet$};
\node at (-5,-2) {$\bullet$};
\node at (-5,-1) {$\bullet$};
\node at (-5,0) {$\bullet$};
\node at (-5,1) {$\circ$};
\node at (-5,2) {$\bullet$};

\node at (-6,-3) {$\bullet$};
\node at (-6,-2) {$\bullet$};
\node at (-6,-1) {$\bullet$};
\node at (-6,0) {$\bullet$};
\node at (-6,1) {$\circ$};
\node at (-6,2) {$\bullet$};

\draw[fill=gray, draw=none, opacity=0.2](0,0)--(-3.3,-3.3)--(-6.6,-3.3)--(-6.6,1)--(0,1)--(0,0);
\end{scope}

\end{tikzpicture}
\caption{Illustration to the proof of Theorem \ref{thm_nonfinal}. Above: dual picture for the first collision in the case when the weight $m$ of the particle sent from the vertex shared by sides $E_1$ and $E_2$ is greater than $1$. We look for the possible position of $\lambda_3,$ the primitive normal direction to the edge $E_3$.  Geometric configuration III of Figure \ref{fig_firstcollision} restricts $\lambda_3$ to belong to the interior of the gray region -- a truncation of the cone with its apex at $\lambda_1$. The condition that positively oriented triangle with vertices $\lambda_2,$ $0$ and $\lambda_3$ is of $A_n$-type restricts $\lambda_3$ to the lattice points depicted by $\circ.$ This leaves a single possibility (marked by $\otimes$) in case of $m>1.$ Observe that the triangle formed by $\lambda_3,$ $0$ and $\lambda_1$ is of area $\frac{1}{2}$, i.e removing (blowing-down) the edge $E_2$ results in a unimodular vertex at the intersection of extensions of sides $E_1$ and $E_3$ of the original polygon. Below: the case of $m=1,$ where there is a single possibility for every natural $n$ -- the mass of the second particle. Again, each of these possibilities is marked by symbol $\otimes$.}\label{fig_nonfinalproof}
\end{figure}

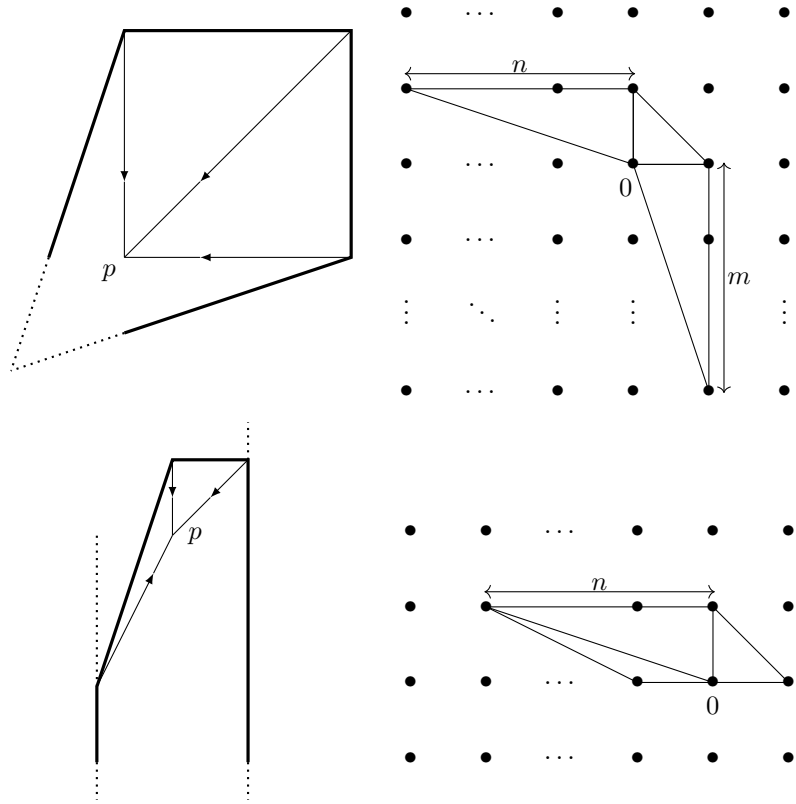
\begin{figure}

\begin{tikzpicture}
\draw[very thick](0,0)--(1,3)--(4,3)--(4,0)--(1,-1);
\draw[thick,dotted](0,0)--(-0.5,-1.5);
\draw[thick,dotted](1,-1)--(-0.5,-1.5);

\draw[-latex](1,3)--(1,1);
\draw(1,1)--(1,0);

\draw[-latex](4,0)--(2,0);
\draw(2,0)--(1,0);

\draw[-latex](4,3)--(2,1);
\draw(2,1)--(1,0);

\node at (0.8,-0.2) {$p$};

\begin{scope}[xshift=220,yshift=35]
\draw(0,0)--(1,0)--(0,1)--(0,0);

\node at (-0.1,-0.3) {$0$};

\draw(0,0)--(0,1)--(-3,1)--(0,0);

\draw[<->](0.025,1.2)--(-3.025,1.2);
\node at (-1.5,1.3) {$n$};

\draw(0,0)--(1,0)--(1,-3)--(0,0);

\draw[<->](1.2,0.025)--(1.2,-3.025);
\node at (1.4,-1.5) {$m$};

\node at (2,-3) {$\bullet$};
\node at (2,-1) {$\bullet$};
\node at (2,0) {$\bullet$};
\node at (2,1) {$\bullet$};
\node at (2,2) {$\bullet$};

\node at (1,-3) {$\bullet$};
\node at (1,-1) {$\bullet$};
\node at (1,0) {$\bullet$};
\node at (1,1) {$\bullet$};
\node at (1,2) {$\bullet$};

\node at (0,-3) {$\bullet$};
\node at (0,-1) {$\bullet$};
\node at (0,0) {$\bullet$};
\node at (0,1) {$\bullet$};
\node at (0,2) {$\bullet$};

\node at (-1,-3) {$\bullet$};
\node at (-1,-1) {$\bullet$};
\node at (-1,0) {$\bullet$};
\node at (-1,1) {$\bullet$};
\node at (-1,2) {$\bullet$};

\node at (-3,-3) {$\bullet$};
\node at (-3,-1) {$\bullet$};
\node at (-3,0) {$\bullet$};
\node at (-3,1) {$\bullet$};
\node at (-3,2) {$\bullet$};

\node at (-2,-3) {$\dots$};
\node at (-2,-1) {$\dots$};
\node at (-2,-0) {$\dots$};
\node at (-2,2) {$\dots$};

\node at (-3,-1.85) {$\vdots$};
\node at (-1,-1.85) {$\vdots$};
\node at (0,-1.85) {$\vdots$};
\node at (2,-1.85) {$\vdots$};

\node at (-2,-1.85) {$\ddots$};
\end{scope}

\begin{scope}[yshift=-190,xshift=75]
\draw[very thick](0,0)--(0,4)--(-1,4)--(-2,1)--(-2,0);
\draw[thick, dotted](0,0)--(0,-0.5);
\draw[thick, dotted](0,4)--(0,4.5);
\draw[thick, dotted](-2,0)--(-2,-0.5);
\draw[thick, dotted](-2,1)--(-2,3);

\draw[-latex](0,4)--(-0.5,3.5);
\draw(-0.5,3.5)--(-1,3);

\draw[-latex](-1,4)--(-1,3.5);
\draw(-1,3.5)--(-1,3);

\draw[-latex](-2,1)--(-1.25,2.5);
\draw(-1.25,2.5)--(-1,3);

\node at (-0.7,3) {$p$};
\end{scope}

\begin{scope}[xshift=250,yshift=-160]
\draw(0,0)--(0,1)--(-3,1)--(0,0);

\draw[<->](0.025,1.2)--(-3.025,1.2);
\node at (-1.5,1.3) {$n$};

\draw(-3,1)--(-1,0)--(0,0);

\draw(0,1)--(1,0)--(0,0);

\node at (0,-0.3) {$0$};

\node at (1,-1) {$\bullet$};
\node at (1,0) {$\bullet$};
\node at (1,1) {$\bullet$};
\node at (1,2) {$\bullet$};

\node at (0,-1) {$\bullet$};
\node at (0,0) {$\bullet$};
\node at (0,1) {$\bullet$};
\node at (0,2) {$\bullet$};

\node at (-1,-1) {$\bullet$};
\node at (-1,0) {$\bullet$};
\node at (-1,1) {$\bullet$};
\node at (-1,2) {$\bullet$};

\node at (-3,-1) {$\bullet$};
\node at (-3,0) {$\bullet$};
\node at (-3,1) {$\bullet$};
\node at (-3,2) {$\bullet$};

\node at (-2,-1) {$\dots$};
\node at (-2,-0) {$\dots$};

\node at (-2,2) {$\dots$};

\node at (-4,-1) {$\bullet$};
\node at (-4,0) {$\bullet$};
\node at (-4,1) {$\bullet$};
\node at (-4,2) {$\bullet$};
\end{scope}

\end{tikzpicture}
\caption{Illustration to the proof of Theorem \ref{thm_nonfinal}. Two possible families of configurations for three consequetive particles, of masses $n,1,m$ and masses $1,n,1,$ and their dual counterparts on the right. Here we use the rigidity result of Figure \ref{fig_nonfinalproof}, i.e. that the middle triangle on the dual picture determines the other two adjacent triangles. In both cases, similarly to I and II of Figure \ref{fig_firstcollision}, we see that at the point of collision $p$ the function $\mathcal{F}_\Phi$ attains its maximum, i.e. the collision is final.}\label{fig_threeparticles}
\end{figure}

\begin{proof}
Consider a vertex $p$ of the tropical caustic of $\Phi,$ corresponding to a collision of particles at a critical time $t_1$ which is less than the final time $t_\Phi.$ By Huygens' principle, we may assume by passing to a positive time $t_1-\varepsilon$ of the original wave front, that the domain $\Phi$ is a polygonal domain, each vertex of which is of $A_n$-type and the vertex $p$ is a first collision (there might be several happening at the same time). There are three geometric configurations as illustrated by Figure \ref{fig_firstcollision}. In cases I and II, we observe that $\mathcal{F}_\Phi$ attains its maximum at $p,$ that is the time of the first collision is final.  The case III is further investigated on Figure \ref{fig_nonfinalproof}, demonstrating certain rigidity, i.e. that the support functions of two adjacent sides together with the weights of the particles determine the support function of the third side. It remains to establish that only two particles can meet at a non-final time. That is illustrated by Figure \ref{fig_threeparticles}, where the above mentioned rigidity phenomenon is used.
\end{proof}

\section{Examples}
\begin{figure}
\includegraphics[width=0.95\textwidth]{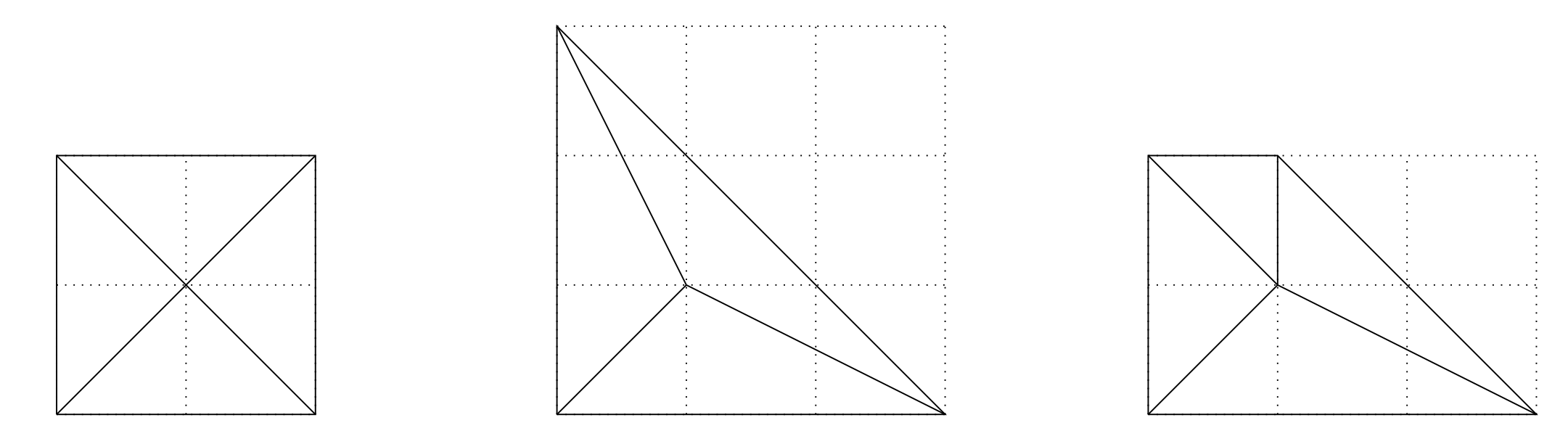}
\caption{Basic examples of polygons with unbranched caustic.}
\label{fig_unbranched}
\end{figure}

\begin{figure}
\includegraphics[width=0.95\textwidth]{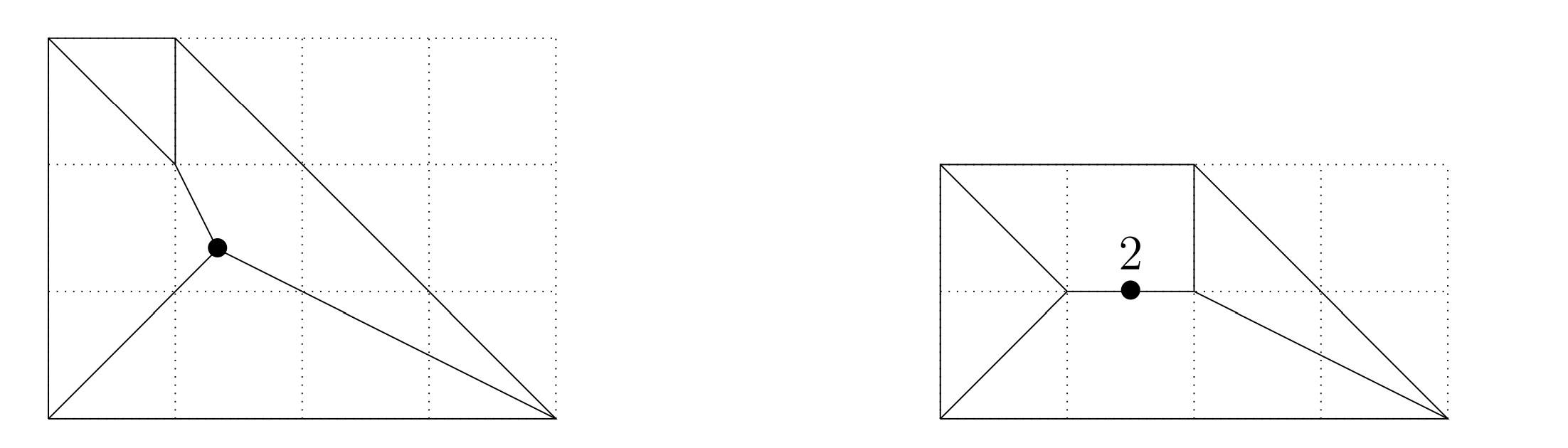}
\caption{The last point of the particle process on the lattice polygon doesn't have to belong to the lattice.}
\label{fig_nonintcenter}
\end{figure}

\begin{figure}
\includegraphics[width=0.95\textwidth]{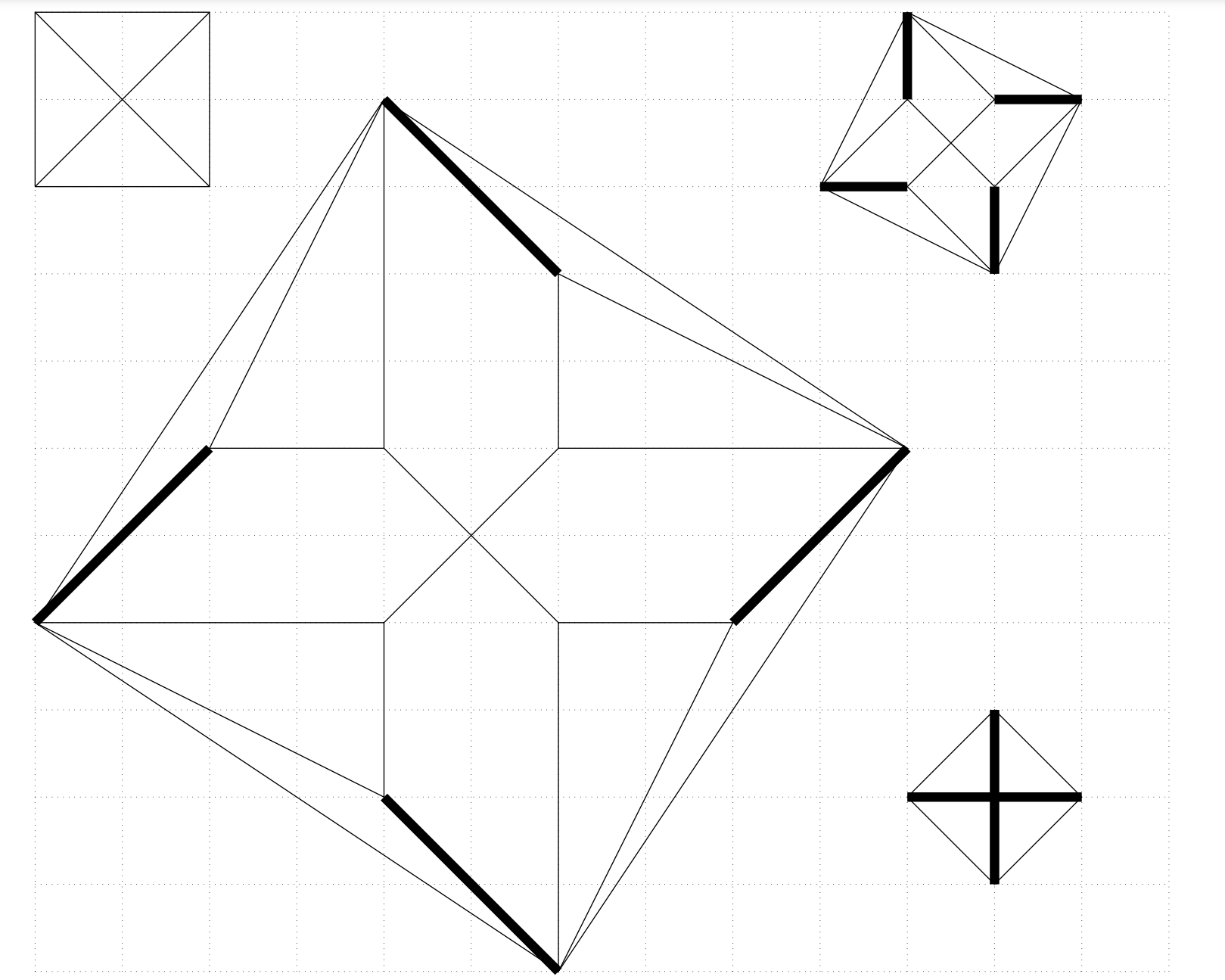}
\caption{A few squares and their caustics. Thick edges are of weight two.}
\label{fig_squares}
\end{figure}

\begin{figure}
\includegraphics[width=0.95\textwidth]{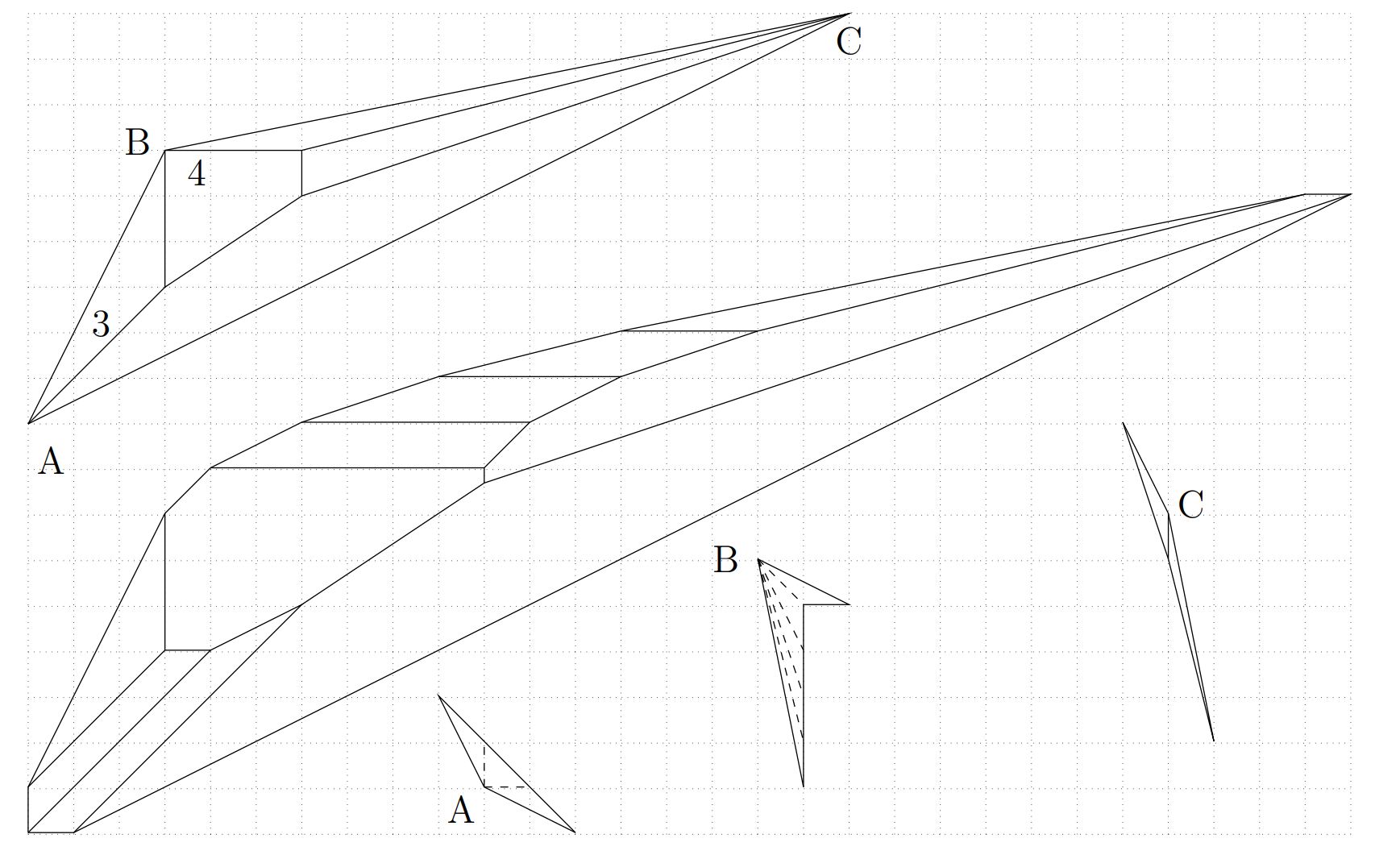}
\caption{Resolution of singularities via caustics.}
\label{fig_resolution}
\end{figure}

The tropical caustic $\mathcal{K}_\Phi$ of a domain $\Phi$ is a finite object, if and only if the domain $\Phi$ is a polygon with rational slope sides. In this case, the easiest method to construct the caustic seems to be to apply the particle process described in the previous section. The case of a unimodular polygon is particularly simple, since we don't need to apply the cone rule of \ref{fig_conerule} -- the only particle sent from every corner of the polygon has mass one and moves with the velocity equal to the sum of primitive vectors going along the adjacent sides. Several very simple examples of this kind may be found on Figures \ref{fig_unbranched} and \ref{fig_nonintcenter}. If the polygon is not unimodular, one needs to apply the cone rule -- an example of a triangle is shown on Figure \ref{fig_resolution}, where the dual pictures to the corners  A, B, and C are shown down on the right. This figure also demonstrates how the caustic can be used to resolve singularities, i.e. to make a unimodular polygon out of a non-unimodular one in a natural way. Next, goes a family of examples where $\Phi$ is a square, see Figure \ref{fig_squares} -- these cases are also simpler, since dual cones to the corners of a square may be identified with the original cones. 

In general, one can use an approximation to the wave front (via iteration of taking the convex hull of the lattice points in the interior) to visualize the caustic, as it is done on Figure \ref{fig_ellipseapprox}. However, this is not how Figure \ref{fig_causticellips} was obtained. Such pictures are built via finding the symplectic minimal model and computing the critical times. By the symplectic minimal model of $\Phi$ we mean such a domain $\hat \Phi$ for which the final time is the same, and its caustic has no vertices outside the final locus, and the caustics of $\Phi$ and $\hat \Phi$ coincide near the final locus (or at infinity, for non-compact domains). For compact domains, minimal models always exist, and are obtained by removing intermediate branchings of the original caustic. For non-compact domains with non-linear asymptotes, such as the standard parabola, the concept of the minimal model doesn't seem to make sense.  The first example we consider in detail is also non-compact, but the asymptotes are linear with rational slope, so the minimal model is a cone.

\subsection{Complement to amoeba of a line.}
Consider an algebraic curve $C$ on the complex algebraic torus $(\mathbb{C}^*)^2.$ Its image under the map $\operatorname{Log}\colon(\mathbb{C}^*)^2\rightarrow(\mathbb{R}^*)^2,$ taking coordinate-wise logarithm, is called the amoeba \cite{GKZ} of $C.$ Each connected component in the complement to the amoeba is convex. Moreover, there is a canonical injective map (see \cite{FPT}) from the set of these connected components to the set of lattice points of the Newton polygon of $C,$ and this construction is $SL_2\mathbb{Z}$ invariant. Therefore, it is natural to study tropical caustics of such connected components and it is reasonable to anticipate that it is possible to recover the Newton polygon together with a lattice point from such a caustic.

\begin{figure}
\includegraphics[width=0.95\textwidth]{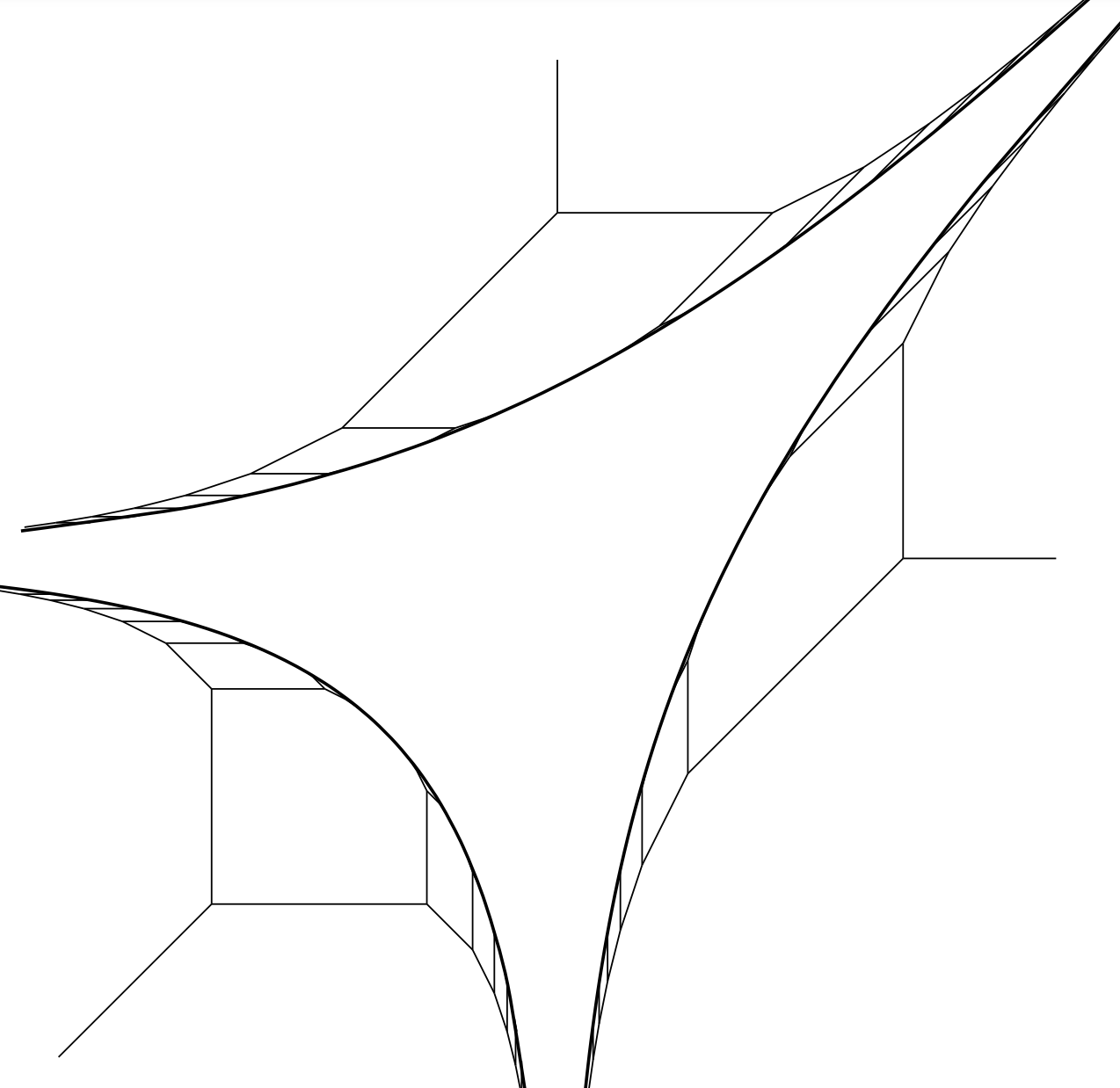}
\caption{Caustics of connected components in the complement to the amoeba of a line.}
\label{fig_amoeba}
\end{figure}

On Figure \ref{fig_amoeba}, the case of an amoeba of a line is presented. We may assume that the line is given by the equation $u+v+1=0,$ where $(u,v)\in(\mathbb{C}^*)^2.$ Note that this amoeba is symmetric under the group with three elements generated by a matrix $\begin{pmatrix}0&-1\\1&1\end{pmatrix}\in SL_2\mathbb{Z},$ thus the three caustics of the connected components in the complement are congruent. We explain how the caustic of the down-left component is constructed. 

First of all, looking from far away, i.e. in the limit under homothety, this component $\Phi$ becomes the quadrant, which corresponds to $\mathbb{C}^2$ from the moment map perspective. The quadrant is the minimal model of the domain under consideration, it caustic is the bisector $\{(s,s):s\leq 0\}.$ If we perform a symplectic blow-up of $\mathbb{C}^2$ at the origin, this corresponds to the corner cut of the quadrant resulting in the polygonal domain $Q_1$ and branching of the caustic. We need to determine the size of such a cut that the domain $\Phi$ is inscribed in $Q_1,$ and therefore the position of the branching of the caustic. If $x,y$ are the coordinates on $\mathbb{R}^2,$ the two support half-planes of $Q$ and $\Phi$ are $-x\geq 0$ and $-y\geq 0.$ If the corner cut is done by the line $-x-y=c,$ then its size is exactly $c$ and the vertex of the caustic of $Q_1$ is located at $(-c,-c).$ To find $c,$ we parametrize the boundary of the domain $\Phi$ as $x<0\mapsto(x,f(x)),$ where $f(x)=\log(1-exp(x))$ and look for $x_1=(f')^{-1}(-1),$ which gives $x_1=-\log(2),$ and therefore $c=-2x_1=2\log(2).$ We proceed further in the similar way, by finding the free coefficients in all tropical monomials contributing to the tropical series $\mathcal{F}_\Phi.$  This is again done by solving the equation $f'(x_{p\over q})=-{p\over q},$ for all rational ${p\over q}$ and finding $c_{p\over q}=-qx_{p\over q}+p \log(1-x_{p\over q}).$  Then, the critical times are equal to $c_{{p_1+p_2}\over {q_1+q_2}}-c_{p_1\over q_1}-c_{p_2\over q_2},$ where $p_1q_2-p_2q_1=1,$ and correspond to the corner cuts by a line with slope $-{p_1+p_2\over q_1+q_2}.$

\subsection{Ellipses}
Consider an ellipse given by $x^2+(\alpha^{-1} y)^2=1,$ for $\alpha>0.$ Then, its support functions (or monomials in the corresponding tropical series) have the form $$px+qy+\sqrt{p^2+(\alpha q)^2}.$$ Therefore, the critical times are computed as $$\sqrt{p_1^2+(\alpha q_1)^2}+\sqrt{p_2^2+(\alpha q_2)^2}-\sqrt{(p_1+p_2)^2+(\alpha (q_1+q_2))^2}.$$  The exact relation between the length's of edges of the caustic and critical times, which is explicitly stated in Remark \ref{rem_crittimes_and_lengths}, is that the lattice length of an edge connecting to vertices is exactly the difference between the times of the two events (particle creation or annihilation) at the ends of the edge. Therefore, to build the caustic it is left to determine the minimal model. For a general ellipse of the form  the corresponding minimal model must have the symmetries $x\mapsto{-x}$ and $y\mapsto{-y}$. There are only two such minimal models, one being a rectangular, and other being a diamond, which is its dual. The latter, however, would have non-final edges of weight two, which is impossible for a domain with a smooth boundary. 

\subsection{Nodal cubic}
\begin{figure}
\includegraphics[width=0.95\textwidth]{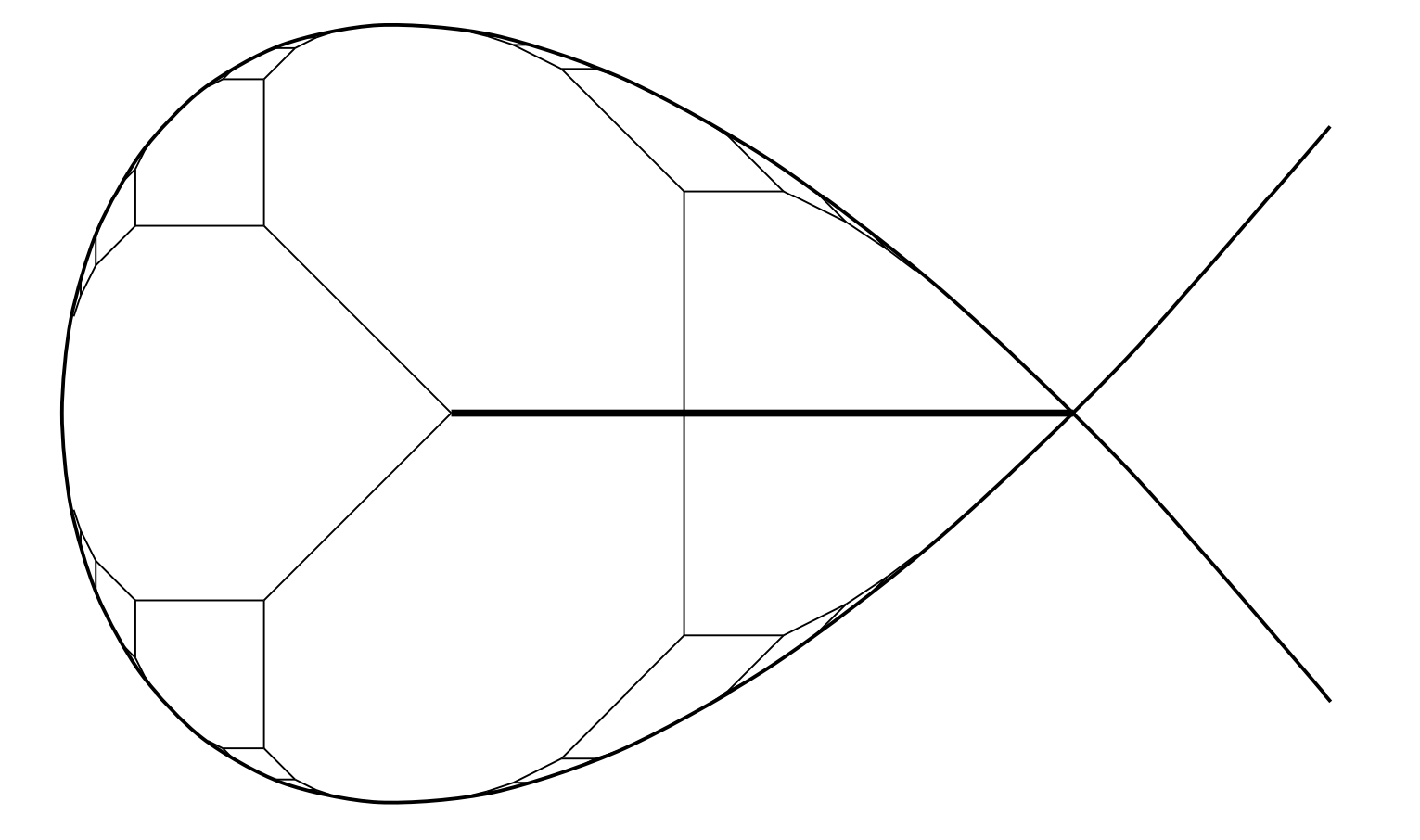}
\caption{A caustic of a convex component in the complement to a real nodal cubic given by equation $y^2=x^2(x+1)$.}
\label{fig_cubic}
\end{figure}

A nodal cubic shown on Figure \ref{fig_cubic} is given by the equation $y^2=x^2(x+1).$ Denote by $\Phi$ the closure of the bounded convex connected component in its complement. The boundary of $\Phi$ is smooth away from the origin. At the origin it is tangent to the bisectors of the upper left and lower left quadrants. In particular, the caustic of $\Phi$ near the origin must consist of a single horizontal edge of weight two. This, and the symmetry $y\mapsto -y$ put constraints on the minimal model of $\Phi$ -- it must be a pentagon of the type shown in the middle of Figure \ref{fig_reflexivepoly} on the right. To find the critical times, and thus the positions of branchings and the lengths of edges, we may use the rational parametrization of the cubic. In the standard way, this parametrization is done by the slopes of lines passing through the origin -- for slopes different from $\pm 1$ such a line intersects the cubic in exactly one point away from the origin. 

\begin{figure}
\includegraphics[width=0.95\textwidth]{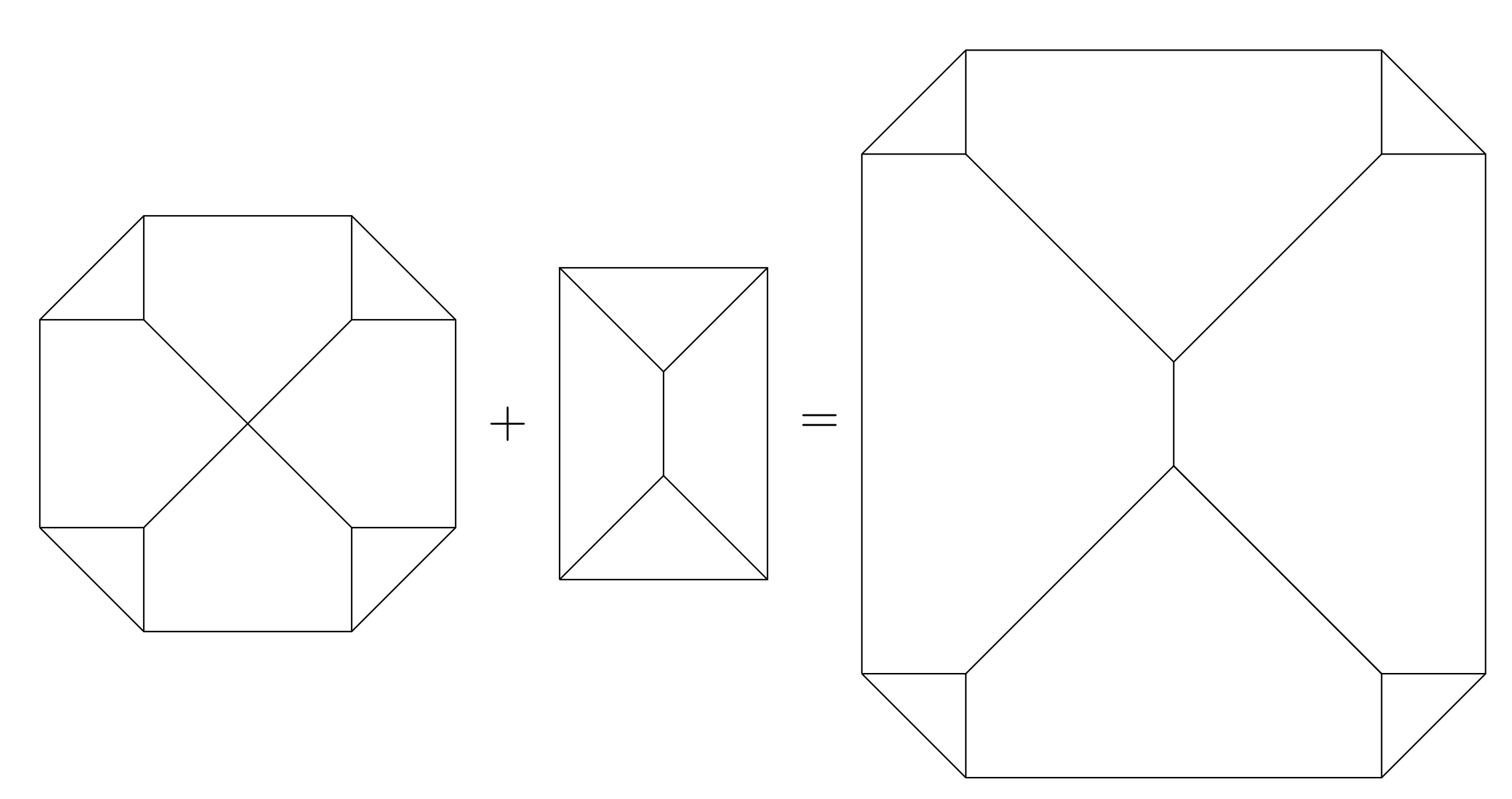}
\caption{Lengths of edges of the caustic are additive when we add to domains with caustics having same final singularity type.}
\label{fig_minksum}
\end{figure}

\begin{figure}
\centering
\includegraphics[width=\textwidth]{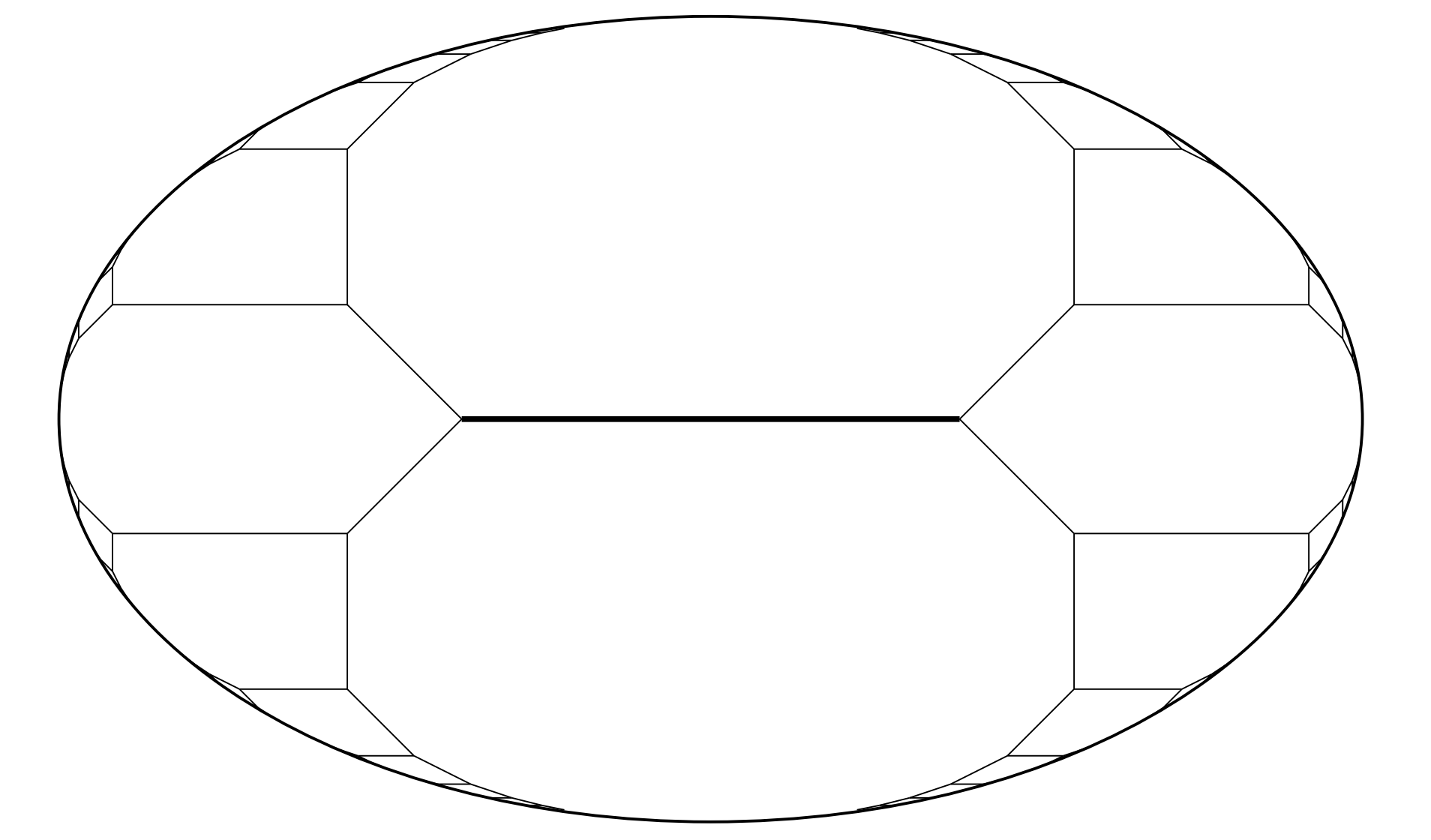}
\caption{A tropical caustic of an ellipse with the aspect $\alpha$ being the golden ratio (see subsection 4.2).}
\label{fig_causticellips}
\end{figure}

\section{Beyond convexity}
In this section, we will briefly describe three strategies to address the question: ``How to extend the notions of tropical wave front and caustic beyond the context of convex domains?''. As before, one would expect that the nature of these geometric objects, i.e. that a wave front is a one-parametric family of curves and a caustic is some weighted graph parametrizing the locus of their special points, is preserved. However, there are two basic properties of caustics that one might need to sacrifice in the higher generality. The first is that a caustic has only positive weights of its edges and the second is that it stays inside the domain. Undoubtedly, the modular invariance of the theory is a virtue to preserve. 

\subsection{Via symplectic geometry of toric surfaces} \label{ssec_toricsurf}
\begin{figure}
\hspace*{-1cm}
\centering
\begin{tikzpicture}

\begin{scope}[scale=0.12]
\node at (-54,0) { };

\draw(-49,-6)--(-49,6)--(-46,12)--(-38,20)--(-35,22)--(-27,26)--(-21,28)--(-17,29)--(-11,30)--(11,30)--(17,29)--(21,28)--(27,26)--(35,22)--(38,20)--(46,12)--(49,6)--(49,-6)--(46,-12)--(38,-20)--(35,-22)--(27,-26)--(21,-28)--(17,-29)--(11,-30)--(-11,-30)--(-17,-29)--(-21,-28)--(-27,-26)--(-35,-22)--(-38,-20)--(-46,-12)--(-49,-6);
\draw(-48,-7)--(-48,7)--(-46,11)--(-36,21)--(-26,26)--(-20,28)--(-16,29)--(16,29)--(20,28)--(26,26)--(36,21)--(46,11)--(48,7)--(48,-7)--(46,-11)--(36,-21)--(26,-26)--(20,-28)--(16,-29)--(-16,-29)--(-20,-28)--(-26,-26)--(-36,-21)--(-46,-11)--(-48,-7);
\draw(-47,-8)--(-47,8)--(-46,10)--(-35,21)--(-25,26)--(-19,28)--(19,28)--(25,26)--(35,21)--(46,10)--(47,8)--(47,-8)--(46,-10)--(35,-21)--(25,-26)--(19,-28)--(-19,-28)--(-25,-26)--(-35,-21)--(-46,-10)--(-47,-8);
\draw(-46,-9)--(-46,9)--(-34,21)--(-24,26)--(-21,27)--(21,27)--(24,26)--(34,21)--(46,9)--(46,-9)--(34,-21)--(24,-26)--(21,-27)--(-21,-27)--(-24,-26)--(-34,-21)--(-46,-9);
\draw(-45,-9)--(-45,9)--(-33,21)--(-23,26)--(23,26)--(33,21)--(45,9)--(45,-9)--(33,-21)--(23,-26)--(-23,-26)--(-33,-21)--(-45,-9);
\draw(-44,-9)--(-44,9)--(-32,21)--(-24,25)--(24,25)--(32,21)--(44,9)--(44,-9)--(32,-21)--(24,-25)--(-24,-25)--(-32,-21)--(-44,-9);
\draw(-43,-9)--(-43,9)--(-31,21)--(-25,24)--(25,24)--(31,21)--(43,9)--(43,-9)--(31,-21)--(25,-24)--(-25,-24)--(-31,-21)--(-43,-9);
\draw(-42,-9)--(-42,9)--(-30,21)--(-26,23)--(26,23)--(30,21)--(42,9)--(42,-9)--(30,-21)--(26,-23)--(-26,-23)--(-30,-21)--(-42,-9);
\draw(-41,-9)--(-41,9)--(-29,21)--(-27,22)--(27,22)--(29,21)--(41,9)--(41,-9)--(29,-21)--(27,-22)--(-27,-22)--(-29,-21)--(-41,-9);
\draw(-40,-9)--(-40,9)--(-28,21)--(28,21)--(40,9)--(40,-9)--(28,-21)--(-28,-21)--(-40,-9);
\draw(-39,-9)--(-39,9)--(-28,20)--(28,20)--(39,9)--(39,-9)--(28,-20)--(-28,-20)--(-39,-9);
\draw(-38,-9)--(-38,9)--(-28,19)--(28,19)--(38,9)--(38,-9)--(28,-19)--(-28,-19)--(-38,-9);
\draw(-37,-9)--(-37,9)--(-28,18)--(28,18)--(37,9)--(37,-9)--(28,-18)--(-28,-18)--(-37,-9);
\draw(-36,-9)--(-36,9)--(-28,17)--(28,17)--(36,9)--(36,-9)--(28,-17)--(-28,-17)--(-36,-9);
\draw(-35,-9)--(-35,9)--(-28,16)--(28,16)--(35,9)--(35,-9)--(28,-16)--(-28,-16)--(-35,-9);
\draw(-34,-9)--(-34,9)--(-28,15)--(28,15)--(34,9)--(34,-9)--(28,-15)--(-28,-15)--(-34,-9);
\draw(-33,-9)--(-33,9)--(-28,14)--(28,14)--(33,9)--(33,-9)--(28,-14)--(-28,-14)--(-33,-9);
\draw(-32,-9)--(-32,9)--(-28,13)--(28,13)--(32,9)--(32,-9)--(28,-13)--(-28,-13)--(-32,-9);
\draw(-31,-9)--(-31,9)--(-28,12)--(28,12)--(31,9)--(31,-9)--(28,-12)--(-28,-12)--(-31,-9);
\draw(-30,-9)--(-30,9)--(-28,11)--(28,11)--(30,9)--(30,-9)--(28,-11)--(-28,-11)--(-30,-9);
\draw(-29,-9)--(-29,9)--(-28,10)--(28,10)--(29,9)--(29,-9)--(28,-10)--(-28,-10)--(-29,-9);
\draw(-28,-9)--(-28,9)--(28,9)--(28,-9)--(-28,-9);
\draw(-27,-8)--(-27,8)--(27,8)--(27,-8)--(-27,-8);
\draw(-26,-7)--(-26,7)--(26,7)--(26,-7)--(-26,-7);
\draw(-25,-6)--(-25,6)--(25,6)--(25,-6)--(-25,-6);
\draw(-24,-5)--(-24,5)--(24,5)--(24,-5)--(-24,-5);
\draw(-23,-4)--(-23,4)--(23,4)--(23,-4)--(-23,-4);
\draw(-22,-3)--(-22,3)--(22,3)--(22,-3)--(-22,-3);
\draw(-21,-2)--(-21,2)--(21,2)--(21,-2)--(-21,-2);
\draw(-20,-1)--(-20,1)--(20,1)--(20,-1)--(-20,-1);
\draw(-19,0)--(19,0)--(-19,0);

\end{scope}
\end{tikzpicture}
\caption{Tropical wave front the ellipse of Figure \ref{fig_causticellips}.}
\label{fig_ellipseapprox}
\end{figure}
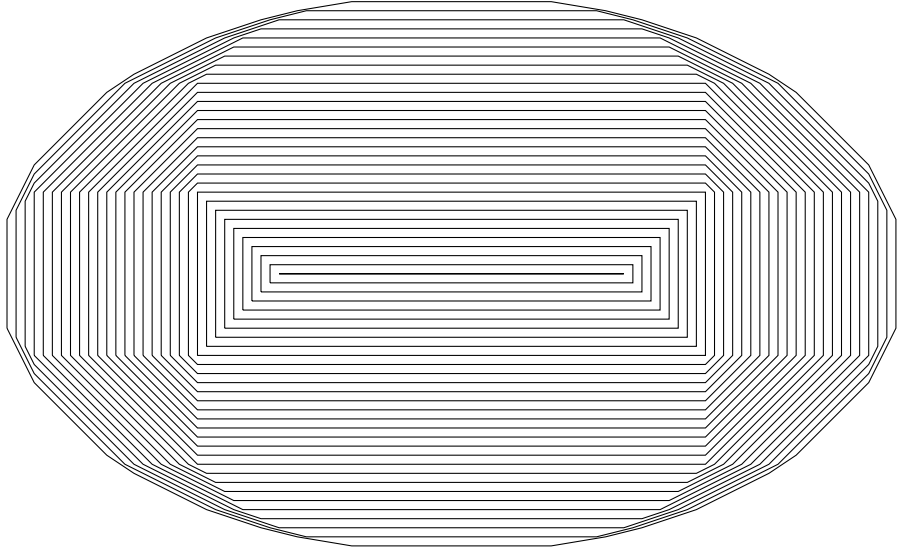

It seems that the most natural of the approaches is through a direct employment of the canonical evolution (\ref{eq_canonicalevolution}) of Theorem \ref{thm_canev} starting with a degenerate form $\omega(0)$. Of course, it would be desirable (but not indispensable) to lift this equation from the de Rham cohomology classes to the actual $2$-forms, which is an interesting problem on its own. I am aware of at least one such generalized situation of mildly degenerate symplectic form for which moment maps are studied, that is of folded symplectic manifolds and, in particular, of origami manifolds \cite{dSGP}, for which the moment domain can fail to be convex and can even be non-simply connected. An obvious limitation of this approach is that by staying within the realm of classical toric varieties one is forced to be content with (possibly non-convex) polygonal rational slope curves as initial states of the tropical wave fronts. There are two ways to overcome this restriction. One is purely analytic in essence and would involve investigating the limits of time-reversed tropical wave front propagation. Another is geometric and would assume developing the theory of moment maps for generalized (such as quantum/non-commutative toric varieties \cite{KLMV21} or pro-orbifold toric surfaces, the construction of which we will describe in a moment equipped with degenerate $2$-forms and proving an analog of Delzant's theorem \cite{D88} for them. 

A pro-orbifold approach to generalized toric surfaces is a straight application of the tropical wave front propagation, that gives a canonical approximation of an arbitrary compact convex domain by rational slope polygons. To be concrete, observe an example illustrated in Figure \ref{fig_ellipseapprox}. Near the final time, a polygon of the tropical propagation is rectangular, it corresponds to the surface $\mathbb{C}P^1\times\mathbb{C}P^1$ with a torus invariant symplectic form. Continuously decreasing the time (and enlarging the polygon) changes this symplectic form, but the surfaces remains the same until the maximal critical time crossing which blows up all four fixed points of the torus action. Continuing going back in time, the surface acquires more and more boundary divisors via blow-ups, and in the limit $t\rightarrow 0$ one gets a pro-manifold model for a symplectic toric surface having an ellipse as its moment domain. In general, as Theorem \ref{thm_An} asserts, surfaces in the approximation may have $A_n$-singularities, therefore the resulting generalized toric surface is rather a pro-orbifold, than pro-manifold.

\subsection{Via lattice approximations} 
\label{ssec_lat_ap}
\begin{figure}
\begin{tikzpicture}
\draw[step=1cm,black,very thin] (0,0.5) grid (4.5,5);
\draw[fill=gray, draw=none, opacity=0.2](0,3)--(2,3)--(2,5)--(5,5)--(5,0)--(0,0);
\draw[fill=gray, draw=none, opacity=0.3](0,2)--(2,2)--(3,3)--(3,5)--(5,5)--(5,0)--(0,0);
\draw[fill=gray, draw=none, opacity=0.3](0,1)--(2,1)--(4,3)--(4,5)--(5,5)--(5,0)--(0,0);
\draw[thick](0,3)--(2,3)--(2,5);
\draw[thick](0,2)--(2,2)--(3,3)--(3,5);
\draw[thick](0,1)--(2,1)--(4,3)--(4,5);

\begin{scope}[xshift=160]
\draw[step=1cm,black,very thin] (0,0.5) grid (6.5,5);
\draw[fill=gray, draw=none, opacity=0.2](0,3)--(2,3)--(4,5)--(7,5)--(7,0)--(0,0);
\draw[fill=gray, draw=none, opacity=0.3](0,2)--(2,2)--(5,5)--(7,5)--(7,0)--(0,0);
\draw[fill=gray, draw=none, opacity=0.3](0,1)--(2,1)--(6,5)--(7,5)--(7,0)--(0,0);
\draw[thick](0,3)--(2,3)--(4,5);
\draw[thick](0,2)--(2,2)--(5,5);
\draw[thick](0,1)--(2,1)--(6,5);
\end{scope}
\end{tikzpicture}
\caption{Limits of lattice approximations (defined in subsection \ref{ssec_lat_ap})  with the standard basis for two $SL(2,\mathbb{Z})$-equivalent non-convex cones. The original cones are shown in light gray, the time $1$ and $2$ propagations are in darker gray.}
\label{fig_lat_ap_ncc}
\end{figure}
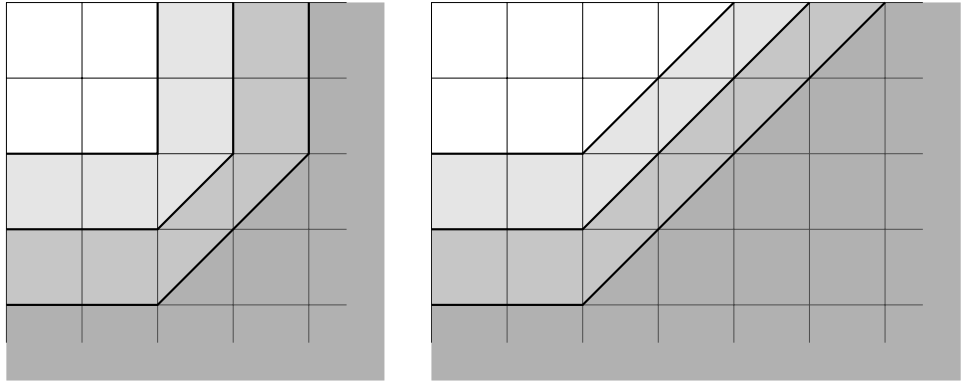

What follows, is an imitation of a different property of a tropical wave front: for a lattice polygon and at the time one it is the boundary of the convex hull of the set of lattice points in the interior of the polygon. For a general domain, there is still a concept of lattice points in the interior. Taking the convex hull, on the other hand, is something to avoid in order to have a possibly non-convex tropical wave front. Another flaw of such a treatment is that it is manifestly discrete in terms of time. Here is a way to remedy this while avoiding the use of the convex hull operation. Denote by $\Phi$ a possibly non-convex domain on $N_\mathbb{R}=N\otimes\mathbb{R},$ where $N$ is a lattice. For a positive $h,$ consider an intersection $\Phi_h$ of $h^{-1}\Phi$ with $N.$ Now, we define moving inside of $\Phi_h$ by $m\in\mathbb{N}$ steps. For this, choose a basis $b$ of $N,$ it turns $N$ into a set of vertices of the corresponding Cayley graph $N_b$ (for a standard square lattice and its standard basis, this gives a planar graph with unit square faces). Define $\Phi_h(b,m)$ as a set of points in $N$ with a distance at least $m$ to the complement of $\Phi_h$ in $N,$ where the distance is computed in terms of the graph distance on $N_b.$ For a fixed $t\geq 0,$ consider a family of discrete sets parametrized by $h>0:$ $$h\Phi_h(b,\floor*{h^{-1}t})\subset\Phi,$$ where $\floor*{-}$ denotes rounding down.  We conjecture that (perhaps under some smoothness assumptions on $\partial\Phi$) this family converges in the Hausdorff sense. Denote this hypothetical limit by $\Phi(b,t).$ It is not very surprising that it depends on the chosen basis $b$ (as Figure \ref{fig_lat_ap_ncc} shows). In order to eliminate this dependence and restore the modular invariance, one can consider an intersection of $\Phi(b,t)$ over all bases $b$ of $N$ as a possible definition for $\Phi(t).$ As a toy, but already not that trivial, example the reader might try working out the case of $\Phi$ being a half-plane (or even a strip) of irrational slope -- the claim is that $\Phi(t)$ is constant with repect to $t$ for such $\Phi$ -- thus extending the notion of tropical wave front to domains which were previously considered non-admissible. 

\subsection{Via extension of caustic coordinates}
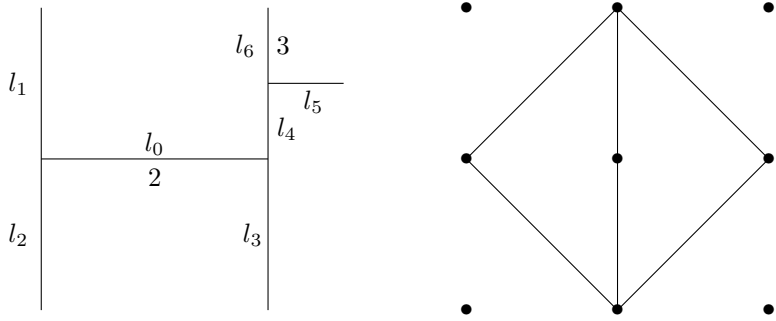
\begin{figure}
\begin{tikzpicture}
\draw(0,0)--(0,4);
\draw(0,2)--(3,2);
\draw(3,0)--(3,4);
\draw(3,3)--(4,3);

\node at (-0.3,1) {$l_2$};
\node at (-0.3,3) {$l_1$};

\node at (1.5,2.2) {$l_0$};
\node at (1.5,1.75) {$2$};

\node at (2.8,1) {$l_3$};
\node at (3.25,2.4) {$l_4$};

\node at (2.7,3.5) {$l_6$};
\node at (3.2,3.5) {$3$};

\node at (3.6,2.75) {$l_5$};

\begin{scope}[xshift=160]
\node at (0,0) {$\bullet$};
\node at (0,2) {$\bullet$};
\node at (0,4) {$\bullet$};
\node at (2,0) {$\bullet$};
\node at (2,2) {$\bullet$};
\node at (2,4) {$\bullet$};
\node at (4,0) {$\bullet$};
\node at (4,2) {$\bullet$};
\node at (4,4) {$\bullet$};

\draw(2,0)--(4,2)--(2,4)--(0,2)--(2,0);
\draw(2,4)--(2,0);
\end{scope}

\end{tikzpicture}
\caption{A metric graph, with lengths of edges denoted by $l_\bullet,$ and decorated by weights (only weights other than one, i.e. two and three, are marked) together with a final singularity type (on the right).}
\label{fig_exampleabstractcasutic}
\end{figure} 

Another approach stems from a closer look at the inversion procedure for the correspondence $\Phi\mapsto\mathcal{K}_\Phi$ used in the proof of Theorem \ref{thm_causticcoordinates}. Let's consider a concrete example. Assume that a caustic has the same final singularity type as the one of an ellipse on Figure \ref{fig_causticellips} and one branching away from the final segment with a weight three edge adjacent to it, Figure \ref{fig_exampleabstractcasutic}. The realizability of such a metric graph as a caustic of some convex domain is equivalent to a set of linear equalities and inequalities on the lengths. The equalities are simply that the distances from the final weight two edge to every endpoint of the graph are equal -- this ensures that a tropical wave front can be reconstructed from the tropical caustic: its vertices at time $t$ are points on the caustic at distance $t$ from the endpoints. Explicitly, it gives $l_1=l_2=l_3=l_4+l_5=l_4+l_6.$

The inequalities are of two types. The first one consists simply of the positivity of each length $l_\bullet$. The second type is non-negativity for lengths of rational slope {\it sides} for the corresponding convex domain. In this example, all but one of them follow from $l_\bullet>0.$ The remaining one is $$l_3+l_4-2l_5\geq 0.$$ 

\begin{figure}
\begin{tikzpicture}

\foreach \x in {0,...,11}
        \foreach \y in {-1,...,8}
	\node at (\x,\y) {$\cdot$};
\draw[fill=gray,thick,opacity=0.2](0,0)--(0,8)--(8,8)--(11,-1)--(11,0)--(0,0);

\draw[dashed](0,0)--(4,4)--(0,8);
\draw[very thick,dashed](4,4)--(7,4);
\draw[dashed](7,4)--(8,5)--(8,8);
\draw[dashed](8,5)--(11,-1);
\draw[dashed](7,4)--(11,0);

\end{tikzpicture}
\caption{A forbidden caustic (shown by dotted lines) and its non-convex domain (shown in light gray filling).}
\label{fig_forbcaustic}
\end{figure}
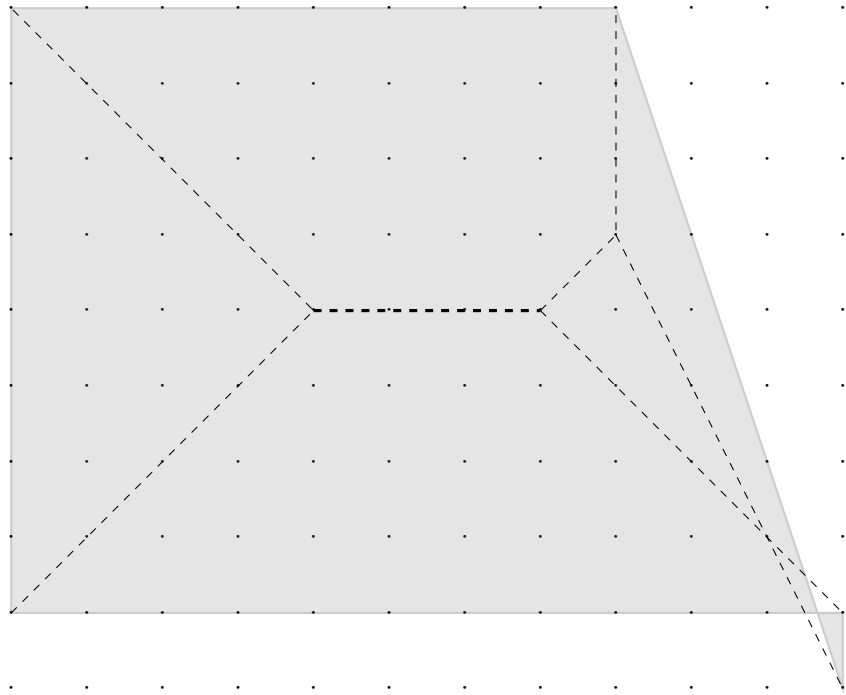

Figure \ref{fig_forbcaustic} shows, what happens if the reversion procedure (for reconstructing a domain from a caustic seen as an abstract metric graph together a singularity type) is applied in the case of this inequality violated. There, the lengths are $l_0=3$, $l_1=l_2=l_3=4,$ $l_4=1$ and $l_5=l_6=3.$ The boundary of the domain is obtained by joining the endpoints of the caustic in their natural cyclic order.

Thus, there is a well-stated problem. What are possible domains for tropical caustics with no inequality restrictions of the second type?

We see that here we already sacrifice one of the basic geometric properties of a tropical caustic that we had in the convex case. Namely, it doesn't necessarily belong to its domain.  Dropping the first set of inequalities for the positivity of edge lengths means allowing a caustic to have negative weights. 

In both situations,  the inversion problem is reduced to some basic linear algebra if caustic curve has only finite number of edges (i.e. if it is semi-algebraic) and to some real hard analysis if the caustic curve is allowed to have infinite number of edges (i.e. if it is analytic).

\subsection{Comparison} 
The three strategies we discussed above are not supposed to be compatible with each other beyond convexity, since already their setups are quite different. Indeed, the second one works with domains as mere subsets of the plane and it is clear that a wave front obtained this way factorizes into a union of independent tropical wave fronts on each of the connected components of the interior of the domain. In contrast to this, the third approach is clearly sensitive to the global structure of the domain as Figure \ref{fig_forbcaustic} demonstrates, where the tropical caustic is {\it not} a union of tropical caustics of the small triangle and the big quadrilateral. The first method is yet another story, in the case of folded symplectic toric manifolds, for example, their moment domains have a finer structure, at least of a multi-set, registering the number of connected components over a given point, and the tropical wave fronts obtained from the canonical evolution (\ref{eq_canonicalevolution}) will definitely reflect it. 

As for the two basic properties of a caustic, that a) it has positive weights and b) belongs to the domain, both are expected to hold in the first approach. Clearly, b) will hold in the second approach, but a) will not. In the third approach,  a) will hold and b) will not, if only the inequalities of the second type are omitted, and both a) and b) will fail if we omit all the inequalities on caustic coordinates.

\section{Concluding remarks}  
There is an obvious question of extending the current framework to higher dimensions. Some work in this direction is already done in \cite{K21}. At the moment of writing, a new article \cite{MS25} about the foundations of tropical wave fronts and caustics in arbitrary dimensions is being prepared. Another matter is about more general ambient spaces. For example, classical caustics are also studied on curved surfaces --- in particular, one of the oldest problems in the field, the last geometric statement of Jacobi, that was resolved in \cite{IK04}, is concerned with the number of cusps of a caustic on an ellipsoid. It would be interesting to examine if problems about classical caustics can be reduced to problems about tropical caustics, as it happens in the relationship between algebraic and tropical geometries. 

Exploring the singularity theory of tropical caustics (especially in the non-convex case) and using it to understand the structure of the space of planar curves and their discriminant is an exciting direction for future research.  Another avenue is applying tropical caustics to algebraic and symplectic geometry. A particular problem is understanding the moduli space of quantum toric varieties \cite{KLMV21}. It is clear already that caustics define a tropical structure on this space. Thus, one could try to compute its tropical (co)homology \cite{IKMZ19}. Perhaps a more feasible question is whether one can extract a canonical ``motivic'' decomposition of convex domains from caustics. We conjecture that by iterating the caustic construction (i.e. taking caustics of the connected components with infinitely many sides in the complements to the original caustic) one obtains a canonical decomposition into hexagons, pentagons, quadrilaterals, and triangles, all with rational slope sides.


\end{document}